\documentclass{article}
\usepackage{amsmath}
\usepackage{amssymb}
\usepackage[curve,arrow]{xy}
\usepackage{theorem}
\theorembodyfont{\itshape}
\newtheorem{theorem}{Theorem}[section]
\newtheorem{corollary}[theorem]{Corollary}
\newtheorem{lemma}[theorem]{Lemma}
\newtheorem{proposition}[theorem]{Proposition}
\newtheorem{conjecture}[theorem]{Conjecture}
\theorembodyfont{\rmfamily\upshape}
\newtheorem{definition}[theorem]{Definition}
\newtheorem{example}[theorem]{Example}
\begin{document}
\newenvironment{proof}[1][,]{\medskip\ifcat,#1
\noindent{\it Proof.\ }\else\noindent{\it Proof of #1.\ }\fi}
{\hfill$\square$\medskip}
\def\e#1\e{\begin{equation}#1\end{equation}}
\def\ea#1\ea{\begin{align}#1\end{align}}
\def\eq#1{{\rm(\ref{#1})}}
\def\dim{\mathop{\rm dim}}
\def\Re{\mathop{\rm Re}}
\def\Im{\mathop{\rm Im}}
\def\rank{\mathop{\rm rank}}
\def\vol{\mathop{\rm vol}}
\def\sind{{\ts\mathop{\text{\rm s-ind}}}}
\def\U{\mathbin{\rm U}}
\def\SU{\mathop{\rm SU}}
\def\sSi{{\smash{\sst\Si}}}
\def\sSii{{\smash{\sst\Si_i}}}
\def\ge{\geqslant} 
\def\le{\leqslant} 
\def\R{\mathbin{\mathbb R}}
\def\Z{\mathbin{\mathbb Z}}
\def\C{\mathbin{\mathbb C}}
\def\CP{\mathbb{CP}}
\def\al{\alpha}
\def\be{\beta}
\def\ga{\gamma}
\def\de{\delta}
\def\ep{\epsilon}
\def\la{\lambda}
\def\th{\theta}
\def\vp{\varphi}
\def\si{\sigma}
\def\om{\omega}
\def\up{\upsilon}
\def\Up{\Upsilon}
\def\De{\Delta}
\def\La{\Lambda}
\def\Om{\Omega}
\def\Ga{\Gamma}
\def\Si{\Sigma}
\def\d{{\rm d}}
\def\pd{\partial}
\def\ts{\textstyle}
\def\sst{\scriptscriptstyle}
\def\w{\wedge}
\def\sm{\setminus}
\def\iy{\infty}
\def\ra{\rightarrow}
\def\longra{\longrightarrow}
\def\t{\times}
\def\ha{{\textstyle\frac{1}{2}}}
\def\ti{\tilde}
\def\bs{\boldsymbol}
\def\ov{\overline}
\def\ovB{\,\overline{\!B}}
\def\ms#1{\vert#1\vert^2}
\def\md#1{\vert #1 \vert}
\def\bmd#1{\big\vert #1 \big\vert}
\def\cnm#1#2{\Vert #1 \Vert_{C^{#2}}} 
\def\lnm#1#2{\Vert #1 \Vert_{L^{#2}}} 
\def\bcnm#1#2{\bigl\Vert #1 \bigr\Vert_{C^{#2}}} 
\title{Singularities of special Lagrangian fibrations \\
and the SYZ Conjecture}
\author{Dominic Joyce \\
Lincoln College, Oxford, OX1 3DR\\
e-mail: {\tt dominic.joyce@lincoln.ox.ac.uk}}
\date{}
\maketitle
\begin{abstract}
The SYZ Conjecture explains Mirror Symmetry between
Calabi--Yau 3-folds $M,\hat M$ in terms of {\it special
Lagrangian fibrations} $f:M\ra B$ and $\hat f:\hat M\ra B$
over the same base $B$, whose fibres are dual 3-tori, except
for singular fibres. This paper studies the {\it singularities}
of special Lagrangian fibrations.

We construct many examples of special Lagrangian fibrations
on open subsets of $\C^3$. The simplest are given explicitly,
and the rest use analytic existence results for
$\U(1)$-invariant special Lagrangian 3-folds in $\C^3$. We
then argue that some features of our examples should also
hold for generic special Lagrangian fibrations of (almost)
Calabi--Yau 3-folds, and draw some conclusions on the SYZ
Conjecture.
\end{abstract}

\section{Introduction}
\label{fi1}

In 1996, Strominger, Yau and Zaslow \cite{SYZ} suggested a
geometrical interpretation of Mirror Symmetry between Calabi--Yau 
3-folds $M,\hat M$ in terms of dual fibrations by special Lagrangian
3-tori, now known as the {\it SYZ Conjecture}. Here is an attempt to 
state it.
\medskip

\noindent{\bf The SYZ Conjecture.} {\it Suppose $M$ and\/ $\hat M$ are 
mirror Calabi--Yau $3$-folds. Then (under some additional conditions) 
there should exist a compact topological\/ $3$-manifold\/ $B$ and 
surjective, continuous maps $f:M\ra B$ and\/ $\hat f:\hat M\ra B$, 
such that
\begin{itemize}
\item[{\rm(i)}] There exists a dense open set\/ $B_0\subset B$, such 
that for each\/ $b\in B_0$, the fibres $f^{-1}(b)$ and\/ $\hat f^{-1}(b)$
are nonsingular special Lagrangian $3$-tori $T^3$ in $M$ and\/ $\hat M$.
Furthermore, $f^{-1}(b)$ and\/ $\hat f^{-1}(b)$ are in some sense
dual to one another.
\item[{\rm(ii)}] For each\/ $b\in\De=B\sm B_0$, the fibres $f^{-1}(b)$ 
and\/ $\hat f^{-1}(b)$ are expected to be singular special Lagrangian
$3$-folds in $M$ and\/~$\hat M$.
\end{itemize}}
\medskip

We call $f$ and $\hat f$ {\it special Lagrangian fibrations}, and
$f^{-1}(b)$, $\hat f^{-1}(b)$ for $b\in\De$ the {\it singular fibres}. 
In this paper we consider the question: what is the nature of the 
`singular fibres' of the fibration, and what do $f,\hat f$ look like 
near the singularities?

The main rigorous results of the paper are the construction and study
in \S\ref{fi5}--\S\ref{fi7} of {\it examples} of special Lagrangian
fibrations on open subsets of $\C^3$. The fibrations of \S\ref{fi5}
are completely explicit, but those in \S\ref{fi6} and \S\ref{fi7}
are constructed using analytic existence results from the author's
series of papers \cite{Joyc2,Joyc3,Joyc4} on special Lagrangian
3-folds in $\C^3$ invariant under the $\U(1)$-action
\begin{equation*}
{\rm e}^{i\th}:(z_1,z_2,z_3)\longmapsto
({\rm e}^{i\th}z_1,{\rm e}^{-i\th}z_2,z_3)
\quad\text{for ${\rm e}^{i\th}\in\U(1)$.}
\end{equation*}

However, the heart of the paper is not the rigorous results but
the discussion and conjecture in \S\ref{fi51}, \S\ref{fi61},
\S\ref{fi74} and \S\ref{fi8}. Here we argue, with justifications
but not full proofs, that various features of our examples should
also be true of special Lagrangian fibrations of (almost) Calabi--Yau
3-folds, especially in the generic case.

In particular, we claim that special Lagrangian fibrations
$f:M\ra B$ will in general not be smooth but only piecewise
smooth, and that the discriminant $\De$ of $f$ is of real
codimension one in $B$ and is typically made up of `ribbons'.
We use this to argue that the version of the SYZ Conjecture
above is too strong, because the discriminants $\De,\hat\De$ of
$f,\hat f$ are not homeomorphic, and so cannot coincide in~$B$.

The paper was originally motivated by the work of Gross
\cite{Gros1,Gros2,Gros3} and Ruan \cite{Ruan1,Ruan2,Ruan3}.
The first version, the preprint math.DG/0011179 in November
2000, consisted mainly of conjectures, and so was not suitable
for publication. It was concerned to refute the widespread
assumption in early papers on the SYZ Conjecture that special
Lagrangian fibrations would be smooth.

In this second version I have used the results of
\cite{Joyc2,Joyc3,Joyc4} to prove many of the conjectures
in the first version. I have also reduced the emphasis on
smoothness of fibrations, as I feel the field has moved on
from two years ago and there is no longer a need to argue
the case.

We begin in \S\ref{fi2} and \S\ref{fi3} by introducing special
Lagrangian geometry and special Lagrangian fibrations. Section
\ref{fi4} reviews the main results of \cite{Joyc2,Joyc3,Joyc4}
on $\U(1)$-invariant special Lagrangian 3-folds in $\C^3$. The
new material is \S\ref{fi5}--\S\ref{fi8}. Section \ref{fi5}
defines two explicit special Lagrangian fibrations $F,F':
\C^3\ra\R^3$ with singular fibres in codimension 1 in~$\R^3$.

Section \ref{fi6} defines a more complicated special Lagrangian
fibration $\hat F:V\ra\R^3$ which models a certain kind of singular
behaviour in codimension 2 in $\R^3$. Section \ref{fi7} constructs
a continuous 1-parameter family of special Lagrangian fibrations
$F^t:V\ra\R^3$ for $t\in[0,1]$, where $F^0$ is smooth, but $F^t$
is not smooth for $t\in(0,1]$. Thus the $F^t$ model how to deform
smooth special Lagrangian fibrations to non-smooth ones.

In each of \S\ref{fi5}--\S\ref{fi7} we also discuss what features
of our examples we expect to hold for special Lagrangian fibrations
of (almost) Calabi--Yau 3-folds, and why. Finally, in \S\ref{fi8}
we explain the picture of smooth special Lagrangian fibrations
built up by Gross and Ruan, and by considering the changes as we
deform from a smooth to a generic fibration, we draw some conclusions
on the SYZ Conjecture.
\medskip

\noindent{\it Acknowledgements.} I would like to thank Mark Gross,
Richard Thomas, Nigel Hitchin and David Morrison for helpful
conversations. I was supported by an EPSRC Advanced Research
Fellowship whilst writing this paper.

\section{Special Lagrangian geometry}
\label{fi2}

We now introduce the idea of special Lagrangian submanifolds
(SL $m$-folds), in two different geometric contexts. First,
in \S\ref{fi21}, we define SL $m$-folds in $\C^m$. Then
\S\ref{fi22} discusses SL $m$-folds in {\it almost Calabi--Yau
$m$-folds}, compact K\"ahler manifolds equipped with a
holomorphic volume form which generalize the idea of
Calabi--Yau manifolds. Finally, section \ref{fi23} considers
the {\it singularities} of SL $m$-folds. The principal references
for this section are Harvey and Lawson \cite{HaLa} and the
author~\cite{Joyc1,Joyc10}.

\subsection{Special Lagrangian submanifolds in $\C^m$}
\label{fi21}

We begin by defining {\it calibrations} and {\it calibrated 
submanifolds}, following Harvey and Lawson~\cite{HaLa}.

\begin{definition} Let $(M,g)$ be a Riemannian manifold. An {\it oriented
tangent $k$-plane} $V$ on $M$ is a vector subspace $V$ of
some tangent space $T_xM$ to $M$ with $\dim V=k$, equipped
with an orientation. If $V$ is an oriented tangent $k$-plane
on $M$ then $g\vert_V$ is a Euclidean metric on $V$, so 
combining $g\vert_V$ with the orientation on $V$ gives a 
natural {\it volume form} $\vol_V$ on $V$, which is a 
$k$-form on~$V$.

Now let $\vp$ be a closed $k$-form on $M$. We say that
$\vp$ is a {\it calibration} on $M$ if for every oriented
$k$-plane $V$ on $M$ we have $\vp\vert_V\le \vol_V$. Here
$\vp\vert_V=\al\cdot\vol_V$ for some $\al\in\R$, and 
$\vp\vert_V\le\vol_V$ if $\al\le 1$. Let $N$ be an 
oriented submanifold of $M$ with dimension $k$. Then 
each tangent space $T_xN$ for $x\in N$ is an oriented
tangent $k$-plane. We say that $N$ is a {\it calibrated 
submanifold} if $\vp\vert_{T_xN}=\vol_{T_xN}$ for all~$x\in N$.
\label{fi2def1}
\end{definition}

It is easy to show that calibrated submanifolds are automatically
{\it minimal submanifolds} \cite[Th.~II.4.2]{HaLa}. Here is the 
definition of special Lagrangian submanifolds in $\C^m$, taken
from~\cite[\S III]{HaLa}.

\begin{definition} Let $\C^m$ have complex coordinates $(z_1,\dots,z_m)$, 
and define a metric $g'$, a real 2-form $\om'$ and a complex $m$-form 
$\Om'$ on $\C^m$ by
\e
\begin{split}
g'=\ms{\d z_1}+\cdots+\ms{\d z_m},\quad
\om'&=\ts\frac{i}{2}(\d z_1\w\d\bar z_1+\cdots+\d z_m\w\d\bar z_m),\\
\text{and}\quad\Om'&=\d z_1\w\cdots\w\d z_m.
\end{split}
\label{fi2eq1}
\e
Then $\Re\Om'$ and $\Im\Om'$ are real $m$-forms on $\C^m$. Let $L$
be an oriented real submanifold of $\C^m$ of real dimension $m$. We
say that $L$ is a {\it special Lagrangian submanifold\/} of $\C^m$,
or {\it SL\/ $m$-fold}\/ for short, if $L$ is calibrated with respect
to $\Re\Om'$, in the sense of Definition~\ref{fi2def1}.
\label{fi2def2}
\end{definition}

Harvey and Lawson \cite[Cor.~III.1.11]{HaLa} give the following
alternative characterization of special Lagrangian submanifolds:

\begin{proposition} Let\/ $L$ be a real $m$-dimensional submanifold 
of\/ $\C^m$. Then $L$ admits an orientation making it into an
SL submanifold of\/ $\C^m$ if and only if\/ $\om'\vert_L\equiv 0$ 
and\/~$\Im\Om'\vert_L\equiv 0$.
\label{fi2prop1}
\end{proposition}

An $m$-dimensional submanifold $L$ in $\C^m$ is called {\it Lagrangian} 
if $\om'\vert_L\equiv 0$. Thus special Lagrangian submanifolds are 
Lagrangian submanifolds satisfying the extra condition that 
$\Im\Om'\vert_L\equiv 0$, which is how they get their name.

\subsection{Almost Calabi--Yau $m$-folds and SL $m$-folds} 
\label{fi22}

We shall define special Lagrangian submanifolds not just in
Calabi--Yau manifolds, as usual, but in the much larger
class of {\it almost Calabi--Yau manifolds}.

\begin{definition} Let $m\ge 2$. An {\it almost Calabi--Yau $m$-fold}, or
{\it ACY\/ $m$-fold}\/ for short, is a quadruple $(M,J,\om,\Om)$ 
such that $(M,J)$ is a compact $m$-dimensional complex manifold,
$\om$ is the K\"ahler form of a K\"ahler metric $g$ on $M$, and
$\Om$ is a non-vanishing holomorphic $(m,0)$-form on~$M$.

We call $(M,J,\om,\Om)$ a {\it Calabi--Yau $m$-fold}, or {\it CY\/ 
$m$-fold}\/ for short, if in addition $\om$ and $\Om$ satisfy
\e
\om^m/m!=(-1)^{m(m-1)/2}(i/2)^m\Om\w\bar\Om.
\label{fi2eq2}
\e
Then for each $x\in M$ there exists an isomorphism $T_xM\cong\C^m$
that identifies $g_x,\om_x$ and $\Om_x$ with the flat versions
$g',\om',\Om'$ on $\C^m$ in \eq{fi2eq1}. Furthermore, $g$ is Ricci-flat
and its holonomy group is a subgroup of~$\SU(m)$.
\label{fi2def3}
\end{definition}

This is not the usual definition of a Calabi--Yau manifold, but
is essentially equivalent to it.

\begin{definition} Let $(M,J,\om,\Om)$ be an almost Calabi--Yau $m$-fold,
and $N$ a real $m$-dimensional submanifold of $M$. We call $N$ a
{\it special Lagrangian submanifold}, or {\it SL $m$-fold} for
short, if $\om\vert_N\equiv\Im\Om\vert_N\equiv 0$. It easily
follows that $\Re\Om\vert_N$ is a nonvanishing $m$-form on $N$.
Thus $N$ is orientable, with a unique orientation in which
$\Re\Om\vert_N$ is positive.
\label{fi2def4}
\end{definition}

Again, this is not the usual definition of special Lagrangian
submanifold, but is essentially equivalent to it. When
$(M,J,\om,\Om)$ is a Calabi--Yau $m$-fold, $N$ is special
Lagrangian if and only if it is calibrated w.r.t.\ $\Re\Om$.
More generally \cite[\S 9.5]{Joyc1}, SL $m$-folds in an ACY
$m$-fold are calibrated w.r.t.\ $\Re\Om$, but for a suitably
{\it conformally rescaled}\/ metric~$g$.

Thus, we could define SL $m$-folds in ACY $m$-folds using
calibrated geometry, as in Definition \ref{fi2def2}. But in
the author's view the definition of SL $m$-folds using the
vanishing of closed forms is more fundamental than that
using calibrated geometry, and so should be taken as the
primary definition.

The {\it deformation theory} of special Lagrangian submanifolds
was studied by McLean \cite[\S 3]{McLe}, who proved the following
result in the Calabi--Yau case. The extension to the ACY case is
described in~\cite[\S 9.5]{Joyc1}.

\begin{theorem} Let\/ $(M,J,\om,\Om)$ be an almost Calabi--Yau $m$-fold,
and\/ $N$ a compact SL\/ $m$-fold in $M$. Then the moduli space
${\cal M}_{\sst N}$ of special Lagrangian deformations of\/ $N$ is a
smooth manifold of dimension $b^1(N)$, the first Betti number of\/~$N$.
\label{fi2thm1}
\end{theorem}

Using similar methods one can prove~\cite[\S 9.3, \S 9.5]{Joyc1}:

\begin{theorem} Let\/ $\bigl\{(M,J_t,\om_t,\Om_t):t\in(-\ep,\ep)\bigr\}$
be a smooth\/ $1$-parameter family of almost Calabi--Yau $m$-folds.
Suppose $N_0$ is a compact SL\/ $m$-fold in $(M,J_0,\om_0,\Om_0)$, with\/
$[\om_t\vert_{N_0}]=0$ in $H^2(N_0,\R)$ and\/ $[\Im\Om_t\vert_{N_0}]=0$
in $H^m(N_0,\R)$ for all\/ $t\in(-\ep,\ep)$. Then $N_0$ extends to a
smooth\/ $1$-parameter family $\bigl\{N_t:t\in(-\de,\de)\bigr\}$,
where $0<\de\le\ep$ and\/ $N_t$ is a compact SL\/ $m$-fold
in~$(M,J_t,\om_t,\Om_t)$.
\label{fi2thm2}
\end{theorem}

\subsection{Singularities of SL $m$-folds}
\label{fi23}

We now summarize some results on {\it singularities} of SL
$m$-folds, taken from the survey \cite{Joyc10} of a series of
papers \cite{Joyc6,Joyc7,Joyc8,Joyc9} by the author. We start
with a definition on {\it SL cones} in $\C^m$, adapted
from~\cite[\S 3.1]{Joyc10}.

\begin{definition} A (singular) SL $m$-fold $C$ in $\C^m$ is called a
{\it cone} if $C=tC$ for all $t>0$, where $tC=\{t\,{\bf x}:{\bf x}
\in C\}$. Let $C$ be a closed SL cone in $\C^m$ for $m>2$ with an
isolated singularity at 0. Then $\Si=C\cap{\cal S}^{2m-1}$ is a
compact, nonsingular $(m\!-\!1)$-submanifold of ${\cal S}^{2m-1}$,
not necessarily connected. Let $g_\sSi$ be the restriction
of $g'$ to $\Si$, where $g'$ is as in~\eq{fi2eq1}.

Let $\De_\sSi$ be the Laplacian on $(\Si,g_\sSi)$. Define
\e
{\mathcal D}_\sSi=\bigl\{\al\in\R:\text{$\al(\al+m-2)$ is
an eigenvalue of $\De_\sSi$}\bigr\}.
\label{fi2eq5}
\e
Then ${\mathcal D}_\sSi$ is a countable, discrete subset of $\R$. Let
$N(\Si)$ be the number of eigenvalues of $\De_\sSi$ in $(0,2m]$,
counted with multiplicity. Let $G$ be the Lie subgroup of $\SU(m)$
preserving $C$. Define the {\it stability index} $\sind(C)$ to be
\e
\sind(C)=N(\Si)-m^2-2m+1+\dim G.
\label{fi2eq6}
\e
Then $\sind(C)\ge 0$. We call $C$ {\it stable} if~$\sind(C)=0$.
\label{fi2def5}
\end{definition}

In \cite[Def.~3.7]{Joyc10} we define {\it conical singularities}
of SL $m$-folds.

\begin{definition} Let $(M,J,\om,\Om)$ be an almost Calabi--Yau $m$-fold
for $m>2$. Suppose $X$ is a compact singular SL $m$-fold in $M$ with
singularities at distinct points $x_1,\ldots,x_n\in X$, and no other
singularities. Fix isomorphisms $\up_i:\C^m\ra T_{x_i}M$ for
$i=1,\ldots,n$ such that $\up_i^*(\om)=\om'$ and $\up_i^*(\Om)=c_i\Om'$,
where $\om',\Om'$ are as in \eq{fi2eq1} and~$c_i>0$.

Let $C_1,\ldots,C_n$ be SL cones in $\C^m$ with isolated singularities
at 0. For $i=1,\ldots,n$ let $\Si_i=C_i\cap{\cal S}^{2m-1}$, and let
$\mu_i\in(2,3)$ with $(2,\mu_i]\cap{\cal D}_\sSii=\emptyset$, where
${\cal D}_\sSii$ is defined in \eq{fi2eq5}. Then we say that $X$ has a
{\it conical singularity} at $x_i$, with {\it rate} $\mu_i$ and
{\it cone} $C_i$ for $i=1,\ldots,n$, if the following holds.

By Darboux' Theorem there exist embeddings $\Up_i:B_R\ra M$
for $i=1,\ldots,n$ satisfying $\Up_i(0)=x_i$, $\d\Up_i\vert_0=\up_i$
and $\Up_i^*(\om)=\om'$, where $B_R$ is the open ball of radius $R$
about 0 in $\C^m$ for some small $R>0$. Define $\iota_i:\Si_i\t(0,R)
\ra B_R$ by $\iota_i(\si,r)=r\si$ for~$i=1,\ldots,n$.

Define $X'=X\sm\{x_1,\ldots,x_n\}$. Then there should exist a
compact subset $K\subset X'$ such that $X'\sm K$ is a union of
open sets $S_1,\ldots,S_n$ with $S_i\subset\Up_i(B_R)$, whose
closures $\bar S_1,\ldots,\bar S_n$ are disjoint in $X$. For
$i=1,\ldots,n$ and some $R'\in(0,R]$ there should exist a smooth
$\phi_i:\Si_i\t(0,R')\ra B_R$ such that $\Up_i\circ\phi_i:\Si_i
\t(0,R')\ra M$ is a diffeomorphism $\Si_i\t(0,R')\ra S_i$, and
\e
\bmd{\nabla^k(\phi_i-\iota_i)}=O(r^{\mu_i-1-k})
\quad\text{as $r\ra 0$ for $k=0,1$.}
\label{fi2eq7}
\e
Here $\nabla$ is the Levi-Civita connection of the cone metric
$\iota_i^*(g')$ on $\Si_i\t(0,R')$, $\md{\,.\,}$ is computed
using $\iota_i^*(g')$. If the cones $C_1,\ldots,C_n$ are
{\it stable} in the sense of Definition \ref{fi2def5}, then
we say that $X$ has {\it stable conical singularities}.
\label{fi2def6}
\end{definition}

In \cite{Joyc7} we study {\it moduli spaces} of SL $m$-folds
with conical singularities. The case when $C_i$ are {\it stable}
is particularly simple,~\cite[Cor.~6.11]{Joyc7}:

\begin{theorem} Suppose $(M,J,\om,\Om)$ is an almost Calabi--Yau
$m$-fold and\/ $X$ a compact SL\/ $m$-fold in $M$ with stable
conical singularities $x_1,\ldots,x_n$. Let\/ ${\cal M}_{\sst X}$ be
the moduli space of deformations of\/ $X$ as an SL\/ $m$-fold with
conical singularities in $M$. Set\/ $X'=X\sm\{x_1,\ldots,x_n\}$,
and let\/ ${\cal I}_{\sst X'}$ be the image of\/ $H^1_{\rm cs}
(X',\R)$ in $H^1(X',\R)$. Then ${\cal M}_{\sst X}$ is a smooth
manifold of dimension~$\dim{\cal I}_{\sst X'}$.
\label{fi2thm3}
\end{theorem}

Here $H^k_{\rm cs}(X',\R),H^k(X',\R)$ are the compactly-supported
and usual de Rham cohomology groups of $X'$. Note the similarity
with Theorem \ref{fi2thm1}. In \cite[Cor.~7.10]{Joyc7}, Theorem
\ref{fi2thm3} is extended to moduli spaces in {\it families} of
almost Calabi--Yau $m$-folds $(M,J^s,\om^s,\Om^s)$. It implies
that SL $m$-folds $X$ with stable conical singularities persist
under small deformations of $(M,J,\om,\Om)$ satisfying some
necessary cohomological conditions.

In \cite{Joyc8,Joyc9} we study {\it desingularizations} of
SL $m$-folds with conical singularities. Here is the basic
idea. Let $M,X$ and $x_i,C_i$ for $i=1,\ldots,n$ be as in
Definition \ref{fi2def6}. Let $L_i$ be an {\it Asymptotically
Conical SL\/ $m$-fold\/} in $\C^m$, asymptotic to $C_i$ at
infinity. Then $tL_i=\{t\,{\bf x}:{\bf x}\in L_i\}$ is also
asymptotic to $C_i$ for all~$t>0$.

We explicitly construct a 1-parameter family of compact,
nonsingular {\it Lagrangian} $m$-folds $N^t$ in $(M,\om)$ for
$t\in(0,\de)$ by gluing $tL_i$ into $X$ at $x_i$, using a
partition of unity. Then we prove using analysis that for
small $t\in(0,\de)$ we can deform $N^t$ to a {\it special\/}
Lagrangian $m$-fold $\smash{\ti N^t}$ in $M$, so that
$\smash{\ti N^t}\ra X$ as $t\ra 0$ in the sense of currents.

The results are complicated, so we will not reproduce them.
Interested readers are advised to consult \cite[\S 7]{Joyc10}.
From them we deduce results on desingularizing SL 3-folds with
$T^2$-cone singularities in \cite[\S 10]{Joyc10}, which will provide
partial proofs of Conjectures \ref{fi5conj} and~\ref{fi6conj}.

\section{Introduction to special Lagrangian fibrations}
\label{fi3}

We begin by defining {\it special Lagrangian fibrations},
following~\cite[Def.~1.4]{Gros2}.

\begin{definition} Let $(M,J,\om,\Om)$ be an almost Calabi--Yau $m$-fold,
and $B$ a Hausdorff topological space. We call $f:M\ra B$ a
{\it special Lagrangian fibration} if $f$ is a continuous,
surjective map, and for all $b\in B$, $f^{-1}(b)$ is the support
of a special Lagrangian integral current $T$ in $M$ with~$\pd T=0$. 

Here {\it integral currents} are meant in the sense of Geometric
Measure Theory. They are a measure-theoretic generalization of
submanifold, including singular submanifolds. Harvey and Lawson
\cite[\S I]{HaLa} frame their discussion of calibrated geometry
in terms of currents, and define calibrated integral currents as
well as calibrated submanifolds. For an introduction to Geometric
Measure Theory, see Morgan~\cite{Morg}.
\label{fi3def1}
\end{definition}

We shall not use much Geometric Measure Theory in this paper.
The point to note is that the fibres $f^{-1}(b)$ are compact
SL $m$-folds in $M$ without boundary, which may have
{\it singularities} of a fairly general kind.

\begin{definition} Let $(M,J,\om,\Om)$ be an almost Calabi--Yau $m$-fold,
and $f:M\ra B$ a special Lagrangian fibration. The {\it fibres}
$N_b$ of $f$ are $N_b=f^{-1}(b)$ for $b\in B$. We call a fibre
$N_b$ {\it nonsingular} if $B$ has the structure of a smooth
real $m$-manifold near $b$, and $f:M\ra B$ is a smooth submersion
along $N_b$. That is, for each $x\in N_b$ the map $\d_xf:T_xM\ra
T_bB$ is surjective. Otherwise we call $N_b$ a {\it singular fibre}.
Define the {\it discriminant} of $f$ to be $\De=\bigl\{b\in B:N_b$
is a singular fibre$\bigr\}$. Roughly speaking, $\De$ is the set of
singular fibres, and $B\sm\De$ the set of nonsingular fibres. It is
easy to show that $B\sm\De$ is open in $B$, and so $\De$ is closed in~$B$.
\label{fi3def2}
\end{definition}

Note that singular fibres $N_b$ may actually sometimes be
{\it nonsingular} submanifolds of $M$. For instance, one can
write down an explicit SL fibration of $T^6/\Z_2$ including
1-parameter families of singular fibres $T^3/\Z_2$, which are
nonsingular as 3-submanifolds, but are `double fibres' of the
fibration.

Using {\it action angle coordinates}, Duistermaat
\cite[Th.~1.1]{Duis} proves:

\begin{proposition} Let\/ $(M,J,\om,\Om)$ be an almost Calabi--Yau
$m$-fold, and\/ $f:M\ra B$ a special Lagrangian fibration.
Then each connected component of a nonsingular fibre $N_b$
is a nonsingular submanifold of\/ $M$ diffeomorphic to~$T^m$.
\label{fi3prop1}
\end{proposition}

The basic idea is that on a nonsingular fibre $N_b$ one
can define a natural action of $T^*_bB=\R^m$, which turns out
to be transitive on connected components. Now arbitrary SL
fibrations $f:M\ra B$ are difficult to study, as we have
little control over their singular behaviour. So it is helpful
to add extra simplifying assumptions. Two such assumptions we
will consider in this paper are that $f$ is {\it smooth}, and
that $f$ is {\it generic}.

\begin{definition} Let $(M,J,\om,\Om)$ be an almost Calabi--Yau $m$-fold,
and $f:M\ra B$ a special Lagrangian fibration. We call $f$ {\it
smooth}\/ if $B$ is a smooth real $m$-manifold, and $f$ a smooth map.
\label{fi3def3}
\end{definition}

Smoothness has strong consequences for the structure of the
discriminant $\De$. From Gross \cite[\S 1]{Gros2} we deduce:

\begin{proposition} Let\/ $M$ be an almost Calabi--Yau $m$-fold, and\/
$f:M\ra B$ a smooth special Lagrangian fibration. Then the
discriminant\/ $\De$ has Hausdorff codimension at least two in~$B$.
\label{fi3prop2}
\end{proposition}

His proof uses the fact that the fibres are both Lagrangian and minimal.
The Lagrangian assumption is used to prove \cite[Prop.~2.2]{Gros2} that
if $x\in N_b$ and $\rank(\d_xf):T_xM\ra T_bB$ is $k$, then $N_b$ contains 
a $k$-dimensional submanifold through $x$ on which $\rank(\d f)$ is $k$. 
But by a result of Almgren, the singularities of a minimal submanifold 
are of Hausdorff codimension at least two. Combining these two shows that
$\rank(\d_xf)$ cannot be $m-1$, so that if $x$ is a singular point of 
$N_b$ then~$\rank(\d_xf)\le m-2$.

Using these ideas, one can show that if $f:M\ra B$ is a smooth
SL fibration of an almost Calabi--Yau 3-fold with discriminant
$\De$, then under good circumstances we expect the following
properties:
\begin{itemize}
\item[(i)] $\De$ is a union $\De_0\cup\De_1$, where $\De_0$ is
a finite set of points, and $\De_1$ a finite set of open intervals.
Essentially, $\De$ is a {\it graph} in~$B$.
\item[(ii)] For each $b\in\De_1$, the singular set of $N_b$ is a
finite number of circles ${\cal S}^1$, and the singularities are
locally modelled on $L\t\R$ in $\C^2\t\C$, where $L$ is a special
Lagrangian 2-fold in $\C^2$ with an isolated singularity at~0.
\end{itemize}
That is, singular fibres occur in {\it codimension two} in the base,
and the generic singular fibre has a {\it one-dimensional} singular
set. Next we define {\it generic} SL fibrations.

\begin{definition} Let $(M,J,\om,\Om)$ be a Calabi--Yau or almost Calabi--Yau 
$m$-fold, and $f:M\ra B$ a special Lagrangian fibration of $(M,J,\om,\Om)$. 
We shall say that some property of $f$ is {\it generic} if for all K\"ahler 
forms $\ti\om$ on $M$ in the same K\"ahler class as $\om$ and sufficiently 
close to $\om$, there exists close to $f$ a special Lagrangian fibration 
$\ti f:M\ra B$ of the almost Calabi--Yau $m$-fold $(M,J,\ti\om,\Om)$ with
the same property. Examples of properties of $f$ that might or might not 
be generic are: existence, smoothness, every singular fibre has only 
finitely many singular points, and so~on.
\label{fi3def4}
\end{definition}

Here is the reasoning behind this definition. We intend to call 
a property of a special Lagrangian fibration {\it generic} if 
it holds for fibrations of all nearby almost Calabi--Yau $m$-folds 
$(M,\ti J,\ti\om,\ti\Om)$. Now if $N$ is a nonsingular fibre of
$f$, then Theorem \ref{fi2thm2} shows that the only obstructions
to finding an SL $m$-fold in $(M,\ti J,\ti\om,\ti\Om)$ near $N$ are 
that~$[\ti\om\vert_N]\equiv[\Im\ti\Om\vert_N]\equiv 0$. 

To ensure this holds, we restrict attention to ACY $m$-folds 
$(M,\ti J,\ti\om,\ti\Om)$ with $[\ti\om]=[\om]$ in $H^2(M,\R)$ and 
$[\Im\ti\Om]=[\Im\Om]$ in $H^3(M,\R)$. But one can show that if 
$[\Im\ti\Om]=[\Im\Om]$ and $(M,J,\Om)$, $(M,\ti J,\ti\Om)$ are close, 
then they are isomorphic. So we may as well fix $\ti J=J$ and 
$\ti\Om=\Om$, and just vary the K\"ahler form $\ti\om$ within the
K\"ahler class of~$\om$.

\section{$\U(1)$-invariant special Lagrangian 3-folds in $\C^3$}
\label{fi4}

We now review the author's three papers \cite{Joyc2,Joyc3,Joyc4}
studying special Lagrangian 3-folds $N$ in $\C^3$ invariant under
the $\U(1)$-action
\begin{equation}                                             
{\rm e}^{i\th}:(z_1,z_2,z_3)\longmapsto
({\rm e}^{i\th}z_1,{\rm e}^{-i\th}z_2,z_3)
\quad\text{for ${\rm e}^{i\th}\in\U(1)$.}
\label{fi4eq1}
\end{equation}
The three papers are briefly surveyed in \cite{Joyc5}. The
results most relevant to this paper are in \cite[\S 8]{Joyc4},
which constructs large families of {\it $\U(1)$-invariant
special Lagrangian fibrations} on open subsets of $\C^3$. These
will be summarized in \S\ref{fi45}, after some introductory
material needed to understand and explain them.

\subsection{Background material from analysis}
\label{fi41}

A closed, bounded, contractible subset $S$ in $\R^n$ will be
called a {\it domain} if the {\it interior} $S^\circ$ of $S$ is
connected with $S=\overline{S^\circ}$, and the {\it boundary}
$\pd S=S\sm S^\circ$ is a compact embedded hypersurface in
$\R^n$. A domain $S$ in $\R^2$ is called {\it strictly convex}
if $S$ is convex and the curvature of $\pd S$ is nonzero at
every point.

Let $S$ be a domain in $\R^n$. Define $C^k(S)$ for $k\ge 0$ to be
the space of continuous functions $f:S\ra\R$ with $k$ continuous
derivatives, and define the norm $\cnm{.}k$ on $C^k(S)$ by
$\cnm{f}k=\sum_{j=0}^k\sup_S\bmd{\pd^jf}$. Then $C^k(S)$ is a Banach
space. Write $C^\iy(S)=\bigcap_{k=0}^\iy C^k(S)$. For $k\ge 0$ and
$\al\in(0,1)$, define the {\it H\"older space} $C^{k,\al}(S)$ to be
the subset of $f\in C^k(S)$ for which
\begin{equation*}
[\pd^k f]_\al=\sup_{x\ne y\in S}
\frac{\bmd{\pd^kf(x)-\pd^kf(y)}}{\md{x-y}^\al}
\end{equation*}
is finite, and define the {\it H\"older norm} on $C^{k,\al}(S)$ 
to be $\cnm{f}{k,\al}=\cnm{f}k+[\pd^kf]_\al$. Again, $C^{k,\al}(S)$
is a Banach space.

A {\it second-order quasilinear operator} $Q:C^2(S)\ra C^0(S)$
is an operator of the form
\e
\bigl(Qu\bigr)(x)=
\sum_{i,j=1}^na^{ij}(x,u,\pd u)\frac{\pd^2u}{\pd x_i\pd x_j}(x)
+b(x,u,\pd u),
\label{fi4eq2}
\e
where $a^{ij}$ and $b$ are continuous maps $S\t\R\t(\R^n)^*\ra\R$,
and $a^{ij}=a^{ji}$ for all $i,j=1,\ldots,n$. We call the functions
$a^{ij}$ and $b$ the {\it coefficients} of $Q$. We call $Q$ 
{\it elliptic} if the symmetric $n\t n$ matrix $(a^{ij})$ is 
positive definite at every point.

A second-order quasilinear operator $Q$ is in {\it divergence form}
if it is written
\begin{equation*}
\bigl(Qu\bigr)(x)=
\sum_{j=1}^n\frac{\pd}{\pd x_j}\bigl(a^j(x,u,\pd u)\bigr)
+b(x,u,\pd u)
\end{equation*}
for functions $a^j\in C^1\bigr(S\t\R\t(\R^n)^*\bigr)$ for
$j=1,\ldots,n$ and $b\in C^0\bigr(S\t\R\t(\R^n)^*\bigr)$.
If $Q$ is in divergence form, we say that integrable functions
$u,f$ are a {\it weak solution} of the equation $Qu=f$ if $u$
is weakly differentiable with weak derivative $\pd u$, and
$a^j(x,u,\pd u),b(x,u,\pd u)$ are integrable with
\begin{equation*}
-\sum_{j=1}^n\int_S\frac{\pd\psi}{\pd x_j}\cdot a^j(x,u,\pd u)
\d{\bf x}+\int_S\psi\cdot b(x,u,\pd u)\d{\bf x}
=\int_S\psi\cdot f\,\d{\bf x}
\end{equation*}
for all $\psi\in C^1(S)$ with~$\psi\vert_{\pd S}\equiv 0$.

If $Q$ is a second-order quasilinear operator, we may interpret
the equation $Qu=f$ in three different senses:
\begin{itemize}
\item We just say that $Qu=f$ if $u\in C^2(S)$, $f\in C^0(S)$
and $Qu=f$ in $C^0(S)$ in the usual way.
\item We say that $Qu=f$ {\it holds with weak derivatives} if $u$
is twice weakly differentiable and $Qu=f$ holds almost everywhere,
defining $Qu$ using weak derivatives.
\item We say that $Qu=f$ {\it holds weakly} if $Q$ is in divergence
form and $u$ is a weak solution of $Qu=f$. Note that this requires
only that $u$ be {\it once} weakly differentiable, and the second
derivatives of $u$ need not exist even weakly.
\end{itemize}

Clearly the first sense implies the second, which implies the third.
If $Q$ is {\it elliptic} and $a^j,b,f$ are suitably regular, one can
usually show that a weak solution to $Qu=f$ is a classical solution,
so that the three senses are equivalent. But for singular equations
that are not elliptic at every point, the three senses are distinct.

\subsection{Finding the equations}
\label{fi42}

Let $N$ be a special Lagrangian 3-fold in $\C^3$ invariant under
the $\U(1)$-action \eq{fi4eq1}. Locally we can write $N$ in the form
\begin{equation}
\begin{split}
N=\bigl\{(z_1,z_2,z_3)\in\C^3:\,& z_1z_2=v(x,y)+iy,\quad z_3=x+iu(x,y),\\
&\ms{z_1}-\ms{z_2}=2a,\quad (x,y)\in S\bigr\},
\end{split}
\label{fi4eq3}
\end{equation}
where $S$ is a domain in $\R^2$, $a\in\R$ and $u,v:S\ra\R$ are
continuous.

Here $\ms{z_1}-\ms{z_2}$ is the {\it moment map} of the $\U(1)$-action
\eq{fi4eq1}, and so $\ms{z_1}-\ms{z_2}$ is constant on any
$\U(1)$-invariant Lagrangian 3-fold in $\C^3$. We choose the constant
to be $2a$. Effectively \eq{fi4eq3} just means that we choose
$x=\Re(z_3)$ and $y=\Im(z_1z_2)$ as local coordinates on the
2-manifold $N/\U(1)$. Then we find~\cite[Prop.~4.1]{Joyc2}:

\begin{proposition} Let\/ $S,a,u,v$ and\/ $N$ be as above. Then
\begin{itemize}
\setlength{\itemsep}{0pt}
\setlength{\parsep}{0pt}
\item[{\rm(a)}] If\/ $a=0$, then $N$ is a (possibly singular)
special Lagrangian $3$-fold in $\C^3$ if\/ $u,v$ are
differentiable and satisfy
\begin{equation}
\frac{\pd u}{\pd x}=\frac{\pd v}{\pd y}
\quad\text{and}\quad
\frac{\pd v}{\pd x}=-2\bigl(v^2+y^2\bigr)^{1/2}\frac{\pd u}{\pd y}\,,
\label{fi4eq4}
\end{equation}
except at points $(x,0)$ in $S$ with\/ $v(x,0)=0$, where $u,v$ 
need not be differentiable. The singular points of\/ $N$ are those
of the form $(0,0,z_3)$, where $z_3=x+iu(x,0)$ for $(x,0)\in S$ 
with\/~$v(x,0)=0$.
\item[{\rm(b)}] If\/ $a\ne 0$, then $N$ is a nonsingular special 
Lagrangian $3$-fold in $\C^3$ if and only if\/ $u,v$ are differentiable
in $S$ and satisfy
\begin{equation}
\frac{\pd u}{\pd x}=\frac{\pd v}{\pd y}\quad\text{and}\quad
\frac{\pd v}{\pd x}=-2\bigl(v^2+y^2+a^2\bigr)^{1/2}\frac{\pd u}{\pd y}\,.
\label{fi4eq5}
\end{equation}
\end{itemize}
\label{fi4prop1}
\end{proposition}

The proof is elementary: at each point ${\bf z}\in N$ we calculate
the tangent space $T_{\bf z}N$ in terms of $\pd u,\pd v$, and use
Proposition \ref{fi2prop1} to find the conditions for $T_{\bf z}N$
to be a special Lagrangian $\R^3$ in $\C^3$. If ${\bf z}=(0,0,z_3)$
then $\d\bigl(\ms{z_1}-\ms{z_2}\bigr)=0$ at $\bf z$, so $\bf z$ is
a singular point of $N$, and $T_{\bf z}N$ does not exist.

Using \eq{fi4eq5} to write $\frac{\pd}{\pd y}\bigl(\frac{\pd u}{\pd x}
\bigr)$ and $\frac{\pd}{\pd x}\bigl(\frac{\pd u}{\pd y}\bigr)$ in
terms of $v$ and setting $\frac{\pd^{\smash{2}}u}{\pd y\pd x}=
\frac{\pd^{\smash{2}}u}{\pd x\pd y}$, we easily
prove~\cite[Prop.~8.1]{Joyc2}:

\begin{proposition} Let\/ $S$ be a domain in $\R^2$ and\/ $u,v\in C^2(S)$
satisfy \eq{fi4eq5} for $a\ne 0$. Then
\begin{equation}
\frac{\pd}{\pd x}\Bigl[\bigl(v^2+y^2+a^2\bigr)^{-1/2}
\frac{\pd v}{\pd x}\Bigr]+2\,\frac{\pd^2v}{\pd y^2}=0.
\label{fi4eq6}
\end{equation}
Conversely, if\/ $v\in C^2(S)$ satisfies \eq{fi4eq6} then 
there exists $u\in C^2(S)$, unique up to addition of a 
constant\/ $u\mapsto u+c$, such that $u,v$ satisfy~\eq{fi4eq5}.
\label{fi4prop2}
\end{proposition}

Now \eq{fi4eq6} is a second order quasilinear elliptic
equation, in divergence form. Thus we can consider
{\it weak solutions} of \eq{fi4eq6} when $a=0$, which need
be only once weakly differentiable. We shall be interested
in solutions of \eq{fi4eq4} with singularities, and the
corresponding SL 3-folds $N$. It will be helpful to define
a class of {\it singular solutions} of~\eq{fi4eq4}.

\begin{definition} Let $S$ be a domain in $\R^2$ and $u,v\in C^0(S)$. We
say that $(u,v)$ is a {\it singular solution} of \eq{fi4eq4} if
\begin{itemize}
\setlength{\itemsep}{0pt}
\setlength{\parsep}{0pt}
\item[(i)] $u,v$ are weakly differentiable, and their weak
derivatives $\frac{\pd u}{\pd x},\frac{\pd u}{\pd y},
\frac{\pd v}{\pd x},\frac{\pd v}{\pd y}$ satisfy~\eq{fi4eq4}.
\item[(ii)] $v$ is a {\it weak solution} of \eq{fi4eq6} with
$a=0$, as in~\S\ref{fi42}.
\item[(iii)] Define the {\it singular points} of $u,v$ to be
the $(x,0)\in S$ with $v(x,0)=0$. Then except at singular points,
$u,v$ are $C^2$ in $S$ and real analytic in $S^\circ$, and
satisfy \eq{fi4eq4} in the classical sense.
\item[(iv)] For $a\in(0,1]$ there exist $u_a,v_a\in C^2(S)$
satisfying \eq{fi4eq5} such that $u_a\ra u$ and $v_a\ra v$
in $C^0(S)$ as~$a\ra 0_+$.
\end{itemize}
\label{fi4def1}
\end{definition}

This list of properties is somewhat arbitrary. The point is that
\cite[\S 8--\S 9]{Joyc3} gives powerful existence and uniqueness
results for solutions $u,v$ of \eq{fi4eq4} satisfying conditions
(i)--(iv) and various boundary conditions on $\pd S$, and all of
(i)--(iv) are useful in different contexts.

\subsection{Examples}
\label{fi43}

The following example is due to Harvey and
Lawson~\cite[\S III.3.A]{HaLa}.

\begin{example} The map $f:\C^3\ra\R^3$ defined by
\e
f:(z_1,z_2,z_3)\longmapsto \bigl(\,\ms{z_1}-\ms{z_2},
\ms{z_1}-\ms{z_3},\Im(z_1z_2z_3)\bigr)
\label{fi4eq7}
\e
is a smooth special Lagrangian fibration of $\C^3$. The fibres of
$f$ are invariant under a subgroup $\U(1)^2$ in $\SU(3)$ acting by
\e
({\rm e}^{i\th_1},{\rm e}^{i\th_2}):(z_1,z_2,z_3)\longmapsto
({\rm e}^{i\th_1}z_1,{\rm e}^{i\th_2}z_2,{\rm e}^{-i(\th_1+\th_2)}z_3),
\label{fi4eq8}
\e
and every $\U(1)^2$-invariant SL 3-fold in $\C^3$ is locally made up
of fibres of $f$.

Calculation shows that the {\it discriminant}\/ of $f$, in the
sense of \S\ref{fi3}, is
\e
\De=\bigl\{(\al,\al,0),(0,-\al,0),(0,0,-\al):\al\ge 0\bigr\}\subset\R^3.
\label{fi4eq9}
\e
It is a {\it trivalent graph}, of codimension two in~$\R^3$.
\label{fi4ex1}
\end{example}

We are interested in a family of particular fibres of $f$ which
decompose into two pieces. Let $a\in\R$, and define
\e
\begin{split}
N_a=\Bigl\{(z_1&,z_2,z_3)\in\C^3:\ms{z_1}-a=\ms{z_2}+a=\ms{z_3}+\md{a},\\
&\Im\bigl(z_1z_2z_3\bigr)=0,\quad \Re\bigl(z_1z_2z_3\bigr)\ge 0\Bigr\}.
\end{split}
\label{fi4eq10}
\e
Then $N_a$ is half of the fibre $f^{-1}(2a,2a,0)$ when $a\ge 0$,
and half of the fibre $f^{-1}(2a,0,0)$ when $a<0$, so $N_a$ is
special Lagrangian.

One can show that $N_a$ is a nonsingular SL 3-fold diffeomorphic
to ${\cal S}^1\t\R^2$ when $a\ne 0$, and $N_0$ is an SL $T^2$-cone
with one singular point at $(0,0,0)$. Note that even though $N_a$
is defined using an inequality $\Re(z_1z_2z_3)\ge 0$, it has no
boundary. This is because the fibres $f^{-1}(2a,2a,0)$ for $a>0$,
and $f^{-1}(2a,0,0)$ for $a<0$, are actually the union of two
nonsingular SL 3-folds ${\cal S}^1\t\R^2$, which intersect in a
circle. The inequality is used to pick out one of these two
SL 3-folds.

By \cite[Th.~5.1]{Joyc2}, these SL 3-folds $N_a$ can be written in
the form~\eq{fi4eq3}.

\begin{theorem} Let\/ $a\in\R$. Then there exist unique
$u_a,v_a:\R^2\ra\R$ such that
\e
\begin{split}
N=\Bigl\{(z_1&,z_2,z_3)\in\C^3:
\Im(z_3)=u_a\bigl(\Re(z_3),\Im(z_1z_2)\bigr),\\
&\Re(z_1z_2)=v_a\bigl(\Re(z_3),\Im(z_1z_2)\bigr),
\quad\ms{z_1}-\ms{z_2}=2a\Bigr\}
\end{split}
\label{fi4eq11}
\e
is the special Lagrangian $3$-fold\/ $N_a$ of\/ \eq{fi4eq10}. Furthermore:
\begin{itemize}
\setlength{\itemsep}{0pt}
\setlength{\parsep}{0pt}
\item[{\rm(a)}] $u_a,v_a$ are smooth on $\R^2$ and satisfy \eq{fi4eq5}, 
except at\/ $(0,0)$ when $a=0$, where they are only continuous.
\item[{\rm(b)}] $u_a(x,y)<0$ when $y>0$ for all\/ $x$, and\/ $u_a(x,0)=0$
for all\/ $x$, and\/ $u_a(x,y)>0$ when $y<0$ for all\/~$x$.
\item[{\rm(c)}] $v_a(x,y)>0$ when $x>0$ for all\/ $y$, and\/ $v_a(0,y)=0$
for all\/ $y$, and\/ $v_a(x,y)<0$ when $x<0$ for all\/~$y$.
\item[{\rm(d)}] $u_a(0,y)=-y\bigl(\md{a}+\sqrt{y^2+a^2}\,\,\bigr)^{-1/2}$ 
for all\/~$y$.
\item[{\rm(e)}] $v_a(x,0)=x\bigl(x^2+2\md{a}\bigr)^{1/2}$ for all\/~$x$.
\item[{\rm(f)}] $u_{-a}\equiv u_a$ and\/ $v_{-a}\equiv v_a$.
\end{itemize}
\label{fi4thm1}
\end{theorem}

In fact \cite[\S 5]{Joyc2} considers only the case $a\ge 0$, but the
case $a<0$ and part (f) follow quickly by exchanging $z_1$ and $z_2$.
Note that although the $N_a$ for $a>0$ and $a<0$ are both diffeomorphic
to ${\cal S}^1\t\R^2$, nonetheless there is a topological change as
$a$ goes from positive to negative, as the fibres undergo a surgery,
a Dehn twist on~${\cal S}^1$.

\subsection{Generating $u,v$ from a potential $f$}
\label{fi44}

In \cite[Prop.~7.1]{Joyc2} we show that solutions $u,v\in C^1(S)$
of \eq{fi4eq5} come from a {\it potential}\/ $f\in C^2(S)$ with
$\frac{\pd f}{\pd y}=u$ and~$\frac{\pd f}{\pd x}=v$.

\begin{proposition} Let\/ $S$ be a domain in $\R^2$ and\/ $u,v\in C^1(S)$
satisfy \eq{fi4eq5} for $a\ne 0$. Then there exists $f\in C^2(S)$
with\/ $\frac{\pd f}{\pd y}=u$, $\frac{\pd f}{\pd x}=v$ and
\begin{equation}
\Bigl(\Bigl(\frac{\pd f}{\pd x}\Bigr)^2+y^2+a^2\Bigr)^{-1/2}
\frac{\pd^2f}{\pd x^2}+2\,\frac{\pd^2f}{\pd y^2}=0.
\label{fi4eq12}
\end{equation}
This $f$ is unique up to addition of a constant, $f\mapsto f+c$.
Conversely, all solutions of\/ \eq{fi4eq12} yield solutions 
of\/~\eq{fi4eq5}. 
\label{fi4prop4}
\end{proposition}

Equation \eq{fi4eq12} is a second-order quasilinear elliptic equation,
singular when $a=0$, which may be written in divergence form. The
following condensation of \cite[Th.~7.6]{Joyc2} and \cite[Th.s 9.20
\& 9.21]{Joyc3} proves existence and uniqueness for the {\it Dirichlet
problem} for~\eq{fi4eq12}.

\begin{theorem} Suppose $S$ is a strictly convex domain in $\R^2$ invariant
under $(x,y)\mapsto(x,-y)$, and\/ $k\ge 0$, $\al\in(0,1)$. Let\/
$a\in\R$ and\/ $\phi\in C^{k+3,\al}(\pd S)$. Then if\/ $a\ne 0$ there
exists a unique $f\in C^{k+3,\al}(S)$ with\/ $f\vert_{\pd S}=\phi$
satisfying \eq{fi4eq12}. If\/ $a=0$ there exists a unique $f\in C^1(S)$
with\/ $f\vert_{\pd S}=\phi$, which is twice weakly differentiable and
satisfies \eq{fi4eq12} with weak derivatives.

Define $u=\frac{\pd f}{\pd y}$ and\/ $v=\frac{\pd f}{\pd x}$. If\/
$a\ne 0$ then $u,v\in C^{k+2,\al}(S)$ satisfy \eq{fi4eq5}, and if\/
$a=0$ then $u,v\in C^0(S)$ are a singular solution of\/ \eq{fi4eq4},
in the sense of Definition \ref{fi4def1}. Furthermore, $f$ depends
continuously in $C^1(S)$, and\/ $u,v$ depend continuously in $C^0(S)$,
on $(\phi,a)$ in~$C^{k+3,\al}(\pd S)\t\R$.
\label{fi4thm2}
\end{theorem}

Combining Proposition \ref{fi4prop1} and Theorem \ref{fi4thm2} gives
existence and uniqueness for a large class of $\U(1)$-invariant SL
3-folds in $\C^3$, with boundary conditions, including {\it singular}
SL 3-folds. It is interesting that this existence and uniqueness is
{\it entirely unaffected} by singularities appearing in~$S^\circ$. 

\subsection{Special Lagrangian fibrations}
\label{fi45}

We can use Theorem \ref{fi4thm2} to construct large families of
{\it special Lagrangian fibrations} of open subsets of $\C^3$
invariant under the $\U(1)$-action \eq{fi4eq1}, including
singular fibres.

\begin{definition} Let $S$ be a strictly convex domain in $\R^2$
invariant under $(x,y)\!\mapsto\!(x,-y)$, let $U$ be an open set
in $\R^3$, and $\al\in(0,1)$. Suppose $\Phi:U\ra C^{3,\al}(\pd S)$
is a continuous map such that if $(a,b,c)\ne(a,b',c')$ in $U$ then
$\Phi(a,b,c)-\Phi(a,b',c')$ has exactly one local maximum and one
local minimum in~$\pd S$.

Let ${\bs\al}=(a,b,c)\in U$, and let $f_{\bs\al}\in
C^{3,\al}(S)$ be the unique (weak) solution of \eq{fi4eq12} with
$f_{\bs\al}\vert_{\pd S}=\Phi({\bs\al})$, which exists by Theorem
\ref{fi4thm2}. Define $u_{\bs\al}=\frac{\pd f_{\bs\al}}{\pd y}$
and $v_{\bs\al}=\frac{\pd f_{\bs\al}}{\pd x}$. Then
$(u_{\bs\al},v_{\bs\al})$ is a solution of \eq{fi4eq5} if
$a\ne 0$, and a singular solution of \eq{fi4eq4} if $a=0$.
Also $u_{\bs\al},v_{\bs\al}$ depend continuously on ${\bs\al}\in U$
in $C^0(S)$, by Theorem~\ref{fi4thm2}.

For each ${\bs\al}=(a,b,c)$ in $U$, define $N_{\bs\al}$ in $\C^3$ by
\e
\begin{split}
N_{\bs\al}=\bigl\{(z_1,z_2,z_3)\in\C^3:\,&
z_1z_2=v_{\bs\al}(x,y)+iy,\quad z_3=x+iu_{\bs\al}(x,y),\\
&\ms{z_1}-\ms{z_2}=2a,\quad (x,y)\in S^\circ\bigr\}.
\end{split}
\label{fi4eq13}
\e
Then $N_{\bs\al}$ is a noncompact SL 3-fold without boundary in $\C^3$,
which is nonsingular if $a\ne 0$, by Proposition~\ref{fi4prop1}.
\label{fi4def2}
\end{definition}

By \cite[Th.~8.2]{Joyc4} the $N_{\bs\al}$ are the fibres of a
{\it special Lagrangian fibration}.

\begin{theorem} In the situation of Definition \ref{fi4def2}, if\/
${\bs\al}\ne{\bs\al}'$ in $U$ then $N_{\bs\al}\cap N_{{\bs\al}'}
=\emptyset$. There exists an open set\/ $V\subset\C^3$ and a
continuous, surjective map $F:V\ra U$ such that\/ $F^{-1}({\bs\al})
=N_{\bs\al}$ for all\/ ${\bs\al}\in U$. Thus, $F$ is a special
Lagrangian fibration of\/ $V\subset\C^3$, which may include
singular fibres.
\label{fi4thm3}
\end{theorem}

The main step in the proof is to show that distinct $N_{\bs\al}$
do not intersect, so that they fibre $V=\bigcup_{{\bs\al}\in U}
N_{\bs\al}$. The tool we use to do this is the following result
\cite[Th.~7.11]{Joyc2}, \cite[Th.~7.10]{Joyc4}:

\begin{theorem} Suppose $S$ is a strictly convex domain in $\R^2$ invariant
under $(x,y)\mapsto(x,-y)$, and\/ $a\in\R$, $k\ge 0$, $\al\in(0,1)$,
and\/ $\phi_1,\phi_2\in C^{k+3,\al}(\pd S)$. Let\/ $u_j,v_j\in C^0(S)$
be the (singular) solution of\/ \eq{fi4eq4} or \eq{fi4eq5} constructed
in Theorem \ref{fi4thm2} from $\phi_j$, for~$j=1,2$.

Suppose $\phi_1-\phi_2$ has $l$ local maxima and\/ $l$
local minima on $\pd S$. Then $(u_1,v_1)-(u_2,v_2)$ has finitely
many zeroes in $S^\circ$, all isolated. Let there be $n$ zeroes
in $S^\circ$ with multiplicities $k_1,\ldots,k_n$.
Then~$\sum_{i=1}^nk_i\le l-1$.
\label{fi4thm4}
\end{theorem}

Here isolated zeroes of $(u_1,v_1)-(u_2,v_2)$ have a
{\it multiplicity}, defined in \cite[Def.~7.1]{Joyc4},
which is a {\it positive integer} by \cite[\S 6.1]{Joyc2}
and \cite[Cor.~7.6]{Joyc4}. The result provides an upper
bound for the number of zeroes of $(u_1,v_2)-(u_2,v_2)$
in $S^\circ$, counted with multiplicity, in terms of
the boundary data~$\phi_1,\phi_2$.

Suppose $\al=(a,b,c)$ and $\al'=(a',b',c')$ are distinct elements
of $U$. If $a\ne a'$ then $N_{\bs\al}\cap N_{{\bs\al}'}=\emptyset$,
since $\ms{z_1}-\ms{z_2}$ is $2a$ on $N_{\bs\al}$ and $2a'$ on
$N_{{\bs\al}'}$. If $a=a'$ then $\Phi({\bs\al})-\Phi({\bs\al}')$
has one local maximum and one local minimum in $\pd S$, by
Definition \ref{fi4def2}. So Theorem \ref{fi4thm4} applies with $l=1$
to show that $(u_{\bs\al},v_{\bs\al})-(u_{{\bs\al}'},v_{{\bs\al}'})$
has no zeroes in $S^\circ$, and again $N_{\bs\al}\cap N_{{\bs\al}'}
=\emptyset$. Thus distinct $N_{\bs\al}$ do not intersect.

For reasons explained in \cite[\S 8]{Joyc4}, we chose to define
$N_{\bs\al}$ in \eq{fi4eq13} over $S^\circ$ rather than $S$, and
so end up with a noncompact SL 3-fold without boundary rather
than a compact SL 3-fold with boundary. The results can be
extended to compact SL 3-folds $\ov N_{\bs\al}$ with boundary,
but it makes the statements rather more complicated, and
introduces new technical problems when $\ov N_{\bs\al}$ has
singularities on its boundary.

There is a simple way \cite[Ex.~8.3]{Joyc4} to produce families
$\Phi$ satisfying Definition \ref{fi4def2}, and thus generate
many SL fibrations of open subsets of~$\C^3$.

\begin{example} Let $S$ be a strictly convex domain in $\R^2$ invariant
under $(x,y)\!\mapsto\!(x,-y)$, let $\al\in(0,1)$ and $\phi\in C^{3,\al}
(\pd S)$. Define $U=\R^3$ and $\Phi:\R^3\ra C^{3,\al}(\pd S)$ by
$\Phi(a,b,c)=\phi+bx+cy$. If $(a,b,c)\ne(a,b',c')$ then
$\Phi(a,b,c)-\Phi(a,b',c')=(b-b')x+(c-c')y\in C^\iy(\pd S)$. As
$b-b',c-c'$ are not both zero and $S$ is strictly convex, it
easily follows that $(b-b')x+(c-c')y$ has one local maximum
and one local minimum in~$\pd S$.

Hence the conditions of Definition \ref{fi4def2} hold for $S,U$
and $\Phi$, and so Theorem \ref{fi4thm3} defines an open set
$V\subset\C^3$ and a special Lagrangian fibration $F:V\ra\C^3$.
One can also show that changing the parameter $c$ in $U=\R^3$
just translates the fibres $N_{\bs\al}$ in $\C^3$, and
$V=\bigl\{(z_1,z_2,z_3)\in\C^3:(\Re z_3,\Im z_1z_2)\in S^\circ\bigr\}$.
\label{fi4ex2}
\end{example}

\subsection{A rough classification of singular points}
\label{fi46}

In \cite[\S 9]{Joyc4} we study {\it singular points} of a
singular solution $u,v$ of~\eq{fi4eq4}.

\begin{definition} Let $S$ be a domain in $\R^2$, and $u,v\in C^0(S)$ a
singular solution of \eq{fi4eq4}, as in Definition \ref{fi4def1}. Suppose
for simplicity that $S$ is invariant under $(x,y)\mapsto(x,-y)$.
Define $u',v'\in C^0(S)$ by $u'(x,y)=u(x,-y)$ and $v'(x,y)=-v(x,-y)$.
Then $u',v'$ is also a singular solution of~\eq{fi4eq4}.

A {\it singular point}, or {\it singularity}, of $(u,v)$ is a point
$(b,0)\in S$ with $v(b,0)=0$. Observe that a singularity of $(u,v)$
is automatically a zero of $(u,v)-(u',v')$. Conversely, a zero of
$(u,v)-(u',v')$ on the $x$-axis is a singularity. A singularity
of $(u,v)$ is called {\it isolated} if it is an isolated zero
of~$(u,v)-(u',v')$.

Define the {\it multiplicity} of an isolated singularity $(b,0)$ of
$(u,v)$ in $S^\circ$ to be the winding number of $(u,v)-(u',v')$
about 0 along the positively oriented circle $\ga_\ep(b,0)$ of
radius $\ep$ about $(b,0)$, where $\ep>0$ is chosen small enough
that $\ga_\ep(b,0)$ lies in $S^\circ$ and $(b,0)$ is the only
zero of $(u,v)-(u',v')$ inside or on $\ga_\ep(b,0)$. By
\cite[Cor.~7.6]{Joyc4}, this multiplicity is a {\it positive
integer}.
\label{fi4def3}
\end{definition}

Under mild conditions the singularities in $S^\circ$ are
isolated, \cite[Th.~9.2]{Joyc4}:

\begin{theorem} Let\/ $S$ be a domain in $\R^2$ invariant under
$(x,y)\mapsto(x,-y)$, and\/ $u,v\in C^0(S)$ a singular solution
of\/ \eq{fi4eq4}. If\/ $u(x,-y)\equiv u(x,y)$ and\/ $v(x,-y)\equiv
-v(x,y)$ then $(u,v)$ is singular along the $x$-axis in $S$, and
the singularities are nonisolated. Otherwise there are at most
countably many singularities of\/ $(u,v)$ in $S^\circ$, all isolated.
\label{fi4thm5}
\end{theorem}

We divide isolated singularities $(b,0)$ into four {\it types},
depending on the behaviour of $v(x,0)$ near~$(b,0)$.

\begin{definition} Let $S$ be a domain in $\R^2$, and $u,v\in C^0(S)$ a
singular solution of \eq{fi4eq4}, as in Definition \ref{fi4def1}.
Suppose $(b,0)$ is an isolated singular point of $(u,v)$ in $S^\circ$.
Then there exists $\ep>0$ such that for $0<\md{x-b}<\ep$ we
have $(x,0)\in S^\circ$ and $v(x,0)\ne 0$. So by continuity $v$
is either positive or negative on each of $(b-\ep,b)\t\{0\}$
and~$(b,b+\ep)\t\{0\}$.
\begin{itemize}
\setlength{\itemsep}{0pt}
\setlength{\parsep}{0pt}
\item[(i)] if $v(x)<0$ for $x\in(b-\ep,b)$ and $v(x)>0$ for
$x\in(b,b+\ep)$ we say the singularity $(b,0)$ is of
{\it increasing type}.
\item[(ii)] if $v(x)>0$ for $x\in(b-\ep,b)$ and $v(x)<0$ for
$x\in(b,b+\ep)$ we say the singularity $(b,0)$ is of
{\it decreasing type}.
\item[(iii)] if $v(x)<0$ for $x\in(b-\ep,b)$ and $v(x)<0$ for
$x\in(b,b+\ep)$ we say the singularity $(b,0)$ is of
{\it maximum type}.
\item[(iv)] if $v(x)>0$ for $x\in(b-\ep,b)$ and $v(x)>0$ for
$x\in(b,b+\ep)$ we say the singularity $(b,0)$ is of
{\it minimum type}.
\end{itemize}
\label{fi4def4}
\end{definition}

The type determines if the multiplicity is even or
odd,~\cite[Prop.~9.4]{Joyc4}.

\begin{proposition} Let\/ $u,v\in C^0(S)$ be a singular solution of\/
\eq{fi4eq4} on a domain $S$ in $\R^2$, and\/ $(b,0)$ be an
isolated singularity of\/ $(u,v)$ in $S^\circ$ with multiplicity
$k$. If\/ $(b,0)$ is of increasing or decreasing type then $k$ is
odd, and if\/ $(b,0)$ is of maximum or minimum type then $k$ is even.
\label{fi4prop5}
\end{proposition}

Theorem \ref{fi4thm4} gives a criterion for finitely many
singularities,~\cite[Th.~9.7]{Joyc4}:\!

\begin{theorem} Suppose $S$ is a strictly convex domain in $\R^2$
invariant under $(x,y)\mapsto(x,-y)$, and\/ $\phi\in C^{k+3,\al}
(\pd S)$ for $k\ge 0$ and\/ $\al\in(0,1)$. Let\/ $u,v$ be the
singular solution of\/ \eq{fi4eq4} in $C^0(S)$ constructed
from $\phi$ in Theorem~$\ref{fi4thm2}$.

Define $\phi'\in C^{k+3,\al}(\pd S)$ by $\phi'(x,y)=-\phi(x,-y)$.
Suppose $\phi-\phi'$ has $l$ local maxima and\/ $l$ local minima
on $\pd S$. Then $(u,v)$ has finitely many singularities in
$S^\circ$. Let there be $n$ singularities in $S^\circ$ with
multiplicities $k_1,\ldots,k_n$. Then~$\sum_{i=1}^nk_i\le l-1$.
\label{fi4thm6}
\end{theorem}

By applying Theorem \ref{fi4thm2} with $S$ the unit disc in $\R^2$
and $\phi$ a linear combination of functions $\sin(j\th),\cos(j\th)$
on the unit circle $\pd S$, we show~\cite[Cor.~10.10]{Joyc4}:

\begin{theorem} There exist examples of singular solutions $u,v$
of\/ \eq{fi4eq4} with isolated singularities of every possible
multiplicity $n\ge 1$, and with both possible types allowed by
Proposition~\ref{fi4prop5}.
\label{fi4thm7}
\end{theorem}

Combining this with Proposition \ref{fi4prop1} gives examples of
SL 3-folds in $\C^3$ with singularities of an {\it infinite number}
of different geometrical/topological types. We also show in
\cite[\S 10.4]{Joyc4} that singular points with multiplicity
$n\ge 1$ occur in {\it real codimension} $n$ in the family of
all SL 3-folds invariant under the $\U(1)$-action \eq{fi4eq1},
in a well-defined sense.

\section{Two model special Lagrangian fibrations}
\label{fi5}

We shall now define two piecewise smooth special Lagrangian
fibrations $F,F':\C^3\ra\R^3$ with singular fibres of codimension
one in $\R^3$. These will be our local models for the most generic
kind of singularity in special Lagrangian fibrations of generic
Calabi--Yau 3-folds. Here is the first.

\begin{theorem} For each\/ $a\in\R$ and\/ $c\in\C$, define $N_{a,c}$
in $\C^3$ by
\e
\begin{split}
N_{a,c}=\Bigl\{(z_1&,z_2,z_3)\in\C^3:\ms{z_1}-a=\ms{z_2}+a=
\ms{z_3-c}+\md{a},\\
&\Im\bigl(z_1z_2(z_3-c)\bigr)=0,\quad
\Re\bigl(z_1z_2(z_3-c)\bigr)\ge 0\Bigr\}.
\end{split}
\label{fi5eq1}
\e
Then $N_{a,c}$ is a nonsingular SL\/ $3$-fold diffeomorphic to
${\cal S}^1\t\R^2$ if\/ $a\ne 0$, and\/ $N_{0,c}$ is an SL\/
$T^2$-cone singular at\/ $(0,0,c)$. Define $F:\C^3\ra\R\t\C$ by
\begin{gather}
F(z_1,z_2,z_3)=(a,b),\quad\text{where}\quad 2a=\ms{z_1}-\ms{z_2}
\label{fi5eq2}\\
\text{and}\quad
b=\begin{cases}
z_3, & a=z_1=z_2=0, \\
z_3-\bar z_1\bar z_2/\md{z_1}, & a\ge 0,\;\> z_1\ne 0, \\
z_3-\bar z_1\bar z_2/\md{z_2}, & a<0.
\end{cases}
\label{fi5eq3}
\end{gather}
Then $F^{-1}(a,c)=N_{a,c}$ for all\/ $a,c\in\R\t\C$, and\/ $F$ is a
continuous, piecewise-smooth SL fibration of\/ $\C^3$, which is not
smooth on~$\md{z_1}=\md{z_2}$.
\label{fi5thm1}
\end{theorem}

\begin{proof} Comparing \eq{fi4eq10} and \eq{fi5eq1} shows that
$N_{a,c}$ is the translation by $(0,0,c)$ of the SL 3-fold $N_a$
of \S\ref{fi43}. Hence $N_{a,c}$ is a nonsingular special Lagrangian
${\cal S}^1\t\R^2$ if $a\ne 0$, and an SL $T^2$-cone singular at
$(0,0,c)$ if $a=0$. It is also easy to see that $F$ is well-defined,
continuous, piecewise smooth, and not smooth on~$\md{z_1}=\md{z_2}$.

One can show from \eq{fi5eq1} that if $(z_1,z_2,z_3)\in N_{a,c}$
then $2a=\ms{z_1}-\ms{z_2}$ and
\begin{equation*}
z_1z_2(z_3-c)=\begin{cases} \md{z_1}\ms{z_2}, & a\ge 0, \\ 
\ms{z_1}\md{z_2}, & a<0. \end{cases}
\end{equation*}
Thus, if $z_1z_2\ne 0$ dividing by $z_1z_2$ and rearranging yields
\begin{equation*}
c=\begin{cases} 
z_3-\md{z_1}\ms{z_2}/(z_1z_2), & a\ge 0, \\
z_3-\ms{z_1}\md{z_2}/(z_1z_2), & a<0.
\end{cases}
\end{equation*}
Using the equations $\ms{z_1}=z_1\bar z_1$ and $\ms{z_2}=z_2\bar z_2$ 
to rewrite these expressions gives the second case of \eq{fi5eq3} when
$z_2\ne 0$ and the third when~$z_1\ne 0$.

If $z_1z_2=0$, equation \eq{fi5eq1} implies that $\ms{z_3-c}=0$, so
$c=z_3$, giving the first case of \eq{fi5eq3}, the second when $z_2=0$
and the third when $z_1=0$. So, if $(z_1,z_2,z_3)\in N_{a,c}$ then we
can recover $a,c$ from $(z_1,z_2,z_3)$ as in \eq{fi5eq2}--\eq{fi5eq3}.
Conversely, for any $(z_1,z_2,z_3)$ in $\C^3$, defining $a,c$ by
\eq{fi5eq2}--\eq{fi5eq3} and reversing the proof above, we find that
$(z_1,z_2,z_3)\in N_{a,c}$. Hence $F^{-1}(a,c)=N_{a,c}$, and $F$ is
a special Lagrangian fibration of~$\C^3$.
\end{proof}

Using Theorem \ref{fi4thm1} we write the fibres $N_{a,c}$ of
$F$ in the form~\eq{fi4eq3}.

\begin{proposition} The SL\/ $3$-folds $N_{a,c}$ of Theorem \ref{fi5thm1}
may be written
\e
\begin{split}
N_{a,c}=\Bigl\{(z_1&,z_2,z_3)\in\C^3:
\Im(z_3)=u_{a,c}\bigl(\Re(z_3),\Im(z_1z_2)\bigr),\\
&\Re(z_1z_2)=v_{a,c}\bigl(\Re(z_3),\Im(z_1z_2)\bigr),
\quad\ms{z_1}-\ms{z_2}=2a\Bigr\},
\end{split}
\label{fi5eq4}
\e
for $u_{a,c},v_{a,c}:\R^2\ra\R$ defined using the functions
$u_a,v_a$ of Theorem \ref{fi4thm1} by
\e
u_{a,c}(x,y)=u_a(x-\Re c,y)+\Im c\quad\text{and}\quad
v_{a,c}(x,y)=v_a(x-\Re c,y).
\label{fi5eq5}
\e
\label{fi5prop1}
\end{proposition}

By applying the involution $(z_1,z_2,z_3)\mapsto(-z_1,z_2,z_3)$ to
$\C^3$ we transform $F$ to a second SL fibration $F'$. The previous
two results quickly yield:

\begin{theorem} For each\/ $a\in\R$ and\/ $c\in\C$, define $N_{a,c}'$
in $\C^3$ by
\e
\begin{split}
N'_{a,c}=\Bigl\{(z_1&,z_2,z_3)\in\C^3:\ms{z_1}-a=
\ms{z_2}+a=\ms{z_3-c}+\md{a},\\
&\Im\bigl(z_1z_2(z_3-c)\bigr)=0,\quad
\Re\bigl(z_1z_2(z_3-c)\bigr)\le 0\Bigr\}.
\end{split}
\label{fi5eq6}
\e
Then $N_{a,c}'$ is a nonsingular SL\/ $3$-fold diffeomorphic to
${\cal S}^1\t\R^2$ if\/ $a\ne 0$, and\/ $N_{0,c}'$ is an SL\/
$T^2$-cone singular at\/ $(0,0,c)$. Define $F':\C^3\ra\R\t\C$ by
\begin{gather}
F'(z_1,z_2,z_3)=(a,b),\quad\text{where}\quad 2a=\ms{z_1}-\ms{z_2}
\label{fi5eq7}\\
\text{and}\quad
b=\begin{cases}
z_3, & a=z_1=z_2=0, \\
z_3+\bar z_1\bar z_2/\md{z_1}, & a\ge 0,\;\> z_1\ne 0, \\
z_3+\bar z_1\bar z_2/\md{z_2}, & a<0.
\end{cases}
\label{fi5eq8}
\end{gather}
Then $(F')^{-1}(a,c)=N_{a,c}'$ for all\/ $a,c\in\R\t\C$, and\/ $F'$ is
a continuous, piecewise-smooth SL fibration of\/ $\C^3$, which is not
smooth on~$\md{z_1}=\md{z_2}$.
\label{fi5thm2}
\end{theorem}

\begin{proposition} The SL\/ $3$-folds $N_{a,c}'$ of Theorem \ref{fi5thm2}
may be written
\e
\begin{split}
N_{a,c}'=\Bigl\{(z_1&,z_2,z_3)\in\C^3:
\Im(z_3)=u_{a,c}'\bigl(\Re(z_3),\Im(z_1z_2)\bigr),\\
&\Re(z_1z_2)=v_{a,c}'\bigl(\Re(z_3),\Im(z_1z_2)\bigr),
\quad\ms{z_1}-\ms{z_2}=2a\Bigr\},
\end{split}
\label{fi5eq9}
\e
for $u_{a,c}',v_{a,c}':\R^2\ra\R$ defined using the functions
$u_a,v_a$ of Theorem \ref{fi4thm1} by
\e
u_{a,c}'(x,y)=-u_a(x-\Re c,y)+\Im c\quad\text{and}\quad
v_{a,c}'(x,y)=-v_a(x-\Re c,y).
\label{fi5eq10}
\e
\label{fi5prop2}
\end{proposition}

\subsection{Discussion}
\label{fi51}

Theorems \ref{fi5thm1} and \ref{fi5thm2} define SL fibrations $F,F'$
of $\C^3$ in which every fibre is invariant under the $\U(1)$-action
\eq{fi4eq1} and is written in the form \eq{fi4eq3}. The singular
fibres of the fibration are $N_{0,c}$ for $c\in\C$, which is singular
only at $(0,0,c)$. Thus the set of singular points of singular fibres
of the fibration is $\bigl\{(0,0,c):c\in\C\bigr\}$, the complex
$z_3$-axis.

However, $F$ is not smooth on the whole real hypersurface $\md{z_1}=
\md{z_2}$, which includes the set of singular points but many other 
points as well. Thus $F$ fails to be smooth not only at singular 
points of singular fibres, but {\it also} at nonsingular points of 
singular fibres. We should understand the non-smoothness of $F$ as
being related not to a singularity at the point in question, but to
a change in the global topology of the whole fibre.

In Theorems \ref{fi5thm1} and \ref{fi5thm2} the base space $B$ is
$\R\t\C$, and the {\it discriminant} $\De$ of Definition \ref{fi3def2}
is $\{0\}\t\C$. Thus, $\De$ is of {\it real codimension one} in $B$.
Now by Proposition \ref{fi3prop2}, for {\it smooth} SL fibrations
$\De$ has {\it Hausdorff codimension two} in $B$. Therefore the
piecewise-smooth SL fibrations of Theorems \ref{fi5thm1} and
\ref{fi5thm2} have very different behaviour to smooth SL fibrations.

We can discuss the singular fibres of $F,F'$ from the point of view
of \S\ref{fi46}. Observe that $(u_0,v_0)$ of Theorem \ref{fi4thm1}
is a singular solution of \eq{fi4eq4} in the sense of Definition
\ref{fi4def1}, with an isolated singularity at $(0,0)$. Using parts
(b) and (c) of Theorem \ref{fi4thm1}, it is easy to show that this
singularity is of {\it multiplicity one} and {\it increasing type},
in the sense of~\S\ref{fi46}.

Now Propositions \ref{fi5prop1} and \ref{fi5prop2} define the
singular solutions $(u_{0,c},v_{0,c})$ and $(u'_{0,c},v'_{0,c})$
for $c\in\C$ in terms of $(u_0,v_0)$. It easily follows that
$(u_{0,c},v_{0,c})$ has an isolated singularity at $(\Re c,0)$
of multiplicity one and increasing type, and $(u_{0,c}',v_{0,c}')$
has an isolated singularity at $(\Re c,0)$ of multiplicity one and
decreasing type. So, the singularities of the fibration $F$ have
{\it multiplicity one} and {\it increasing type}, and those of
$F'$ have {\it multiplicity one} and {\it decreasing type}, in
the sense of~\S\ref{fi46}.

We can also discuss the fibrations using the framework of
\S\ref{fi23}. The singular fibres of $F$ and $F'$ are special
Lagrangian $T^2$-cones modelled on
\e
C=\bigl\{(z_1,z_2,z_3):\md{z_1}\!=\!\md{z_2}\!=\!\md{z_3},\;
\Im(z_1z_2z_3)\!=\!0,\; \Re(z_1z_2z_3)\!\ge\!0\bigr\}.
\label{fi5eq12}
\e
By \cite[Ex.~3.5]{Joyc10}, $C$ is {\it stable} in the sense of
Definition \ref{fi2def5}. Hence, SL 3-folds with conical
singularities with cone $C$ have a very well-behaved deformation
theory, by Theorem \ref{fi2thm3}. We make our first conjecture.

\begin{conjecture} The SL fibrations $F,F'$ of Theorems \ref{fi5thm1},
\ref{fi5thm2} are generic local models for codimension one
singularities of SL fibrations $f:M\ra B$ of almost Calabi--Yau
$3$-folds, using `generic' as in Definition~\ref{fi3def4}.
\label{fi5conj}
\end{conjecture}

Here I am being deliberately vague on what I mean by a {\it local
model\/} for an SL fibration. Roughly speaking, I want $f$ to have
the same topological structure as $F$ or $F'$ locally, for the
fibres of $f$ to approximate those of $F$ or $F'$ near the singular
points, and for the singular fibres to have conical singularities
with cone $C$, in the sense of Definition~\ref{fi2def6}.

We can use the results of \cite{Joyc10} to give a {\it partial
proof\/} of Conjecture \ref{fi5conj}. Let $(M,J,\om,\Om)$ be an
almost Calabi--Yau 3-fold, and $f:M\ra B$ an SL fibration with
generic fibre $T^3$, which is locally modelled on $F$ or $F'$ near
some singular point $x\in M$. Let $b=f(x)$ and $X=f^{-1}(b)$ be the
fibre through $x$. Then $X$ is locally modelled on $C$ near~$x$.

Suppose that $b$ is generic in the discriminant of $F$ or $F'$, and
that $X$ has only finitely many singularities $x=x_1,\ldots,x_n$,
each modelled on $C$. Then \cite[Th.~4.5]{Joyc10} shows that $X$ has
conical singularities at $x_1,\ldots,x_n$ with cone $C$, in the sense
of Definition \ref{fi2def6}. The moduli space ${\cal M}_{\sst X}$ of
Theorem \ref{fi2thm3} locally coincides with the discriminant of $F$
or $F'$. Hence it is smooth with~$\dim{\cal M}_{\sst X}=\dim{\cal
I}_{\sst X'}=2$.

Now let $\ti\om$ be a K\"ahler form on $M$ close to $\om$ and
in the same K\"ahler class, as in Definition \ref{fi3def4}. As
$C$ is {\it stable}, \cite[Cor.~7.10]{Joyc7} shows that for
$\ti\om$ sufficiently close to $\om$ there exists a deformation
$\ti X$ of $X$, which is an SL $m$-fold in $(M,J,\ti\om,\Om)$
with conical singularities $\ti x_1,\ldots,\ti x_n$ and cone
$C$, and lies in a smooth moduli space ${\cal M}_{\sst\ti X}$
with~$\dim{\cal M}_{\sst\ti X}=2$.

We can then use the desingularization results of \cite{Joyc8,Joyc9},
following \cite[\S 10]{Joyc10}, to show that $\ti X$ must admit two
families ${\cal F}_{\sst\ti X}^\pm$ of {\it desingularizations},
which are SL $T^3$'s in $(M,J,\ti\om,\Om)$, corresponding to the
regions $a>0$, $a<0$ in $F$ or $F'$. We expect that ${\cal M}_{\sst
\ti X}$, ${\cal F}_{\sst\ti X}^+$ and ${\cal F}_{\sst\ti X}^-$
will locally form the fibres of an {\it SL fibration} $\ti f:M\ra B$
of $(M,J,\ti\om,\Om)$ close to $f$ and modelled on $F$ or $F'$
near $x$, for $\ti\om$ sufficiently close to~$\om$.

Finally we discuss {\it holomorphic discs} with boundary in the
fibres of~$F,F'$.

\begin{lemma} Let\/ $a>0$ and\/ $c\in\C$. Then $\bigl\{(z_1,0,c):
z_1\in\C$, $\md{z_1}^2\le 2a\bigr\}$ is a holomorphic disc in
$\C^3$ with boundary in both $N_{a,c}$ and\/ $N'_{a,c}$, of
area~$2\pi a$.

Similarly, let\/ $a<0$ and\/ $c\in\C$. Then $\bigl\{(0,z_2,c):
z_2\in\C$, $\md{z_2}^2\le -2a\bigr\}$ is a holomorphic disc in
$\C^3$ with boundary in both $N_{a,c}$ and\/ $N'_{a,c}$, of
area~$-2\pi a$.
\label{fi5lem}
\end{lemma}

The proof is trivial. Now if $D$ is a holomorphic disc in an
almost Calabi--Yau 3-fold $(M,J,\om,\Om)$ with boundary in an SL
3-fold $L$ in $M$, then the area of $D$ is $\int_D\om=[\om]\cdot[D]$,
where $[\om]\in H^2(M,L;\R)$ is the relative de Rham cohomology
class of $\om$, and $[D]\in H_2(M,L;\Z)$ is the relative homology
class of $D$. Thus the area of a holomorphic disc is essentially a
{\it topological invariant}.

Holomorphic discs $D$ with boundary in SL 3-folds $L$ are typically
stable objects which persist under small deformations of $L$.
Also, the area of $D$ is always positive. Suppose we deform $L$
in $M$ so that the area $[\om]\cdot[D]$ of $D$ shrinks to zero.
Then $D$ shrinks down to a point, so that the boundary $\pd D$
in $L$, a circle, is collapsed to a point. So $L$ develops a
singularity, a $T^2$-cone, by collapsing an ${\cal S}^1$ in $L$
to a point.

We can think of the singularities of the fibrations $F,F'$ as
occurring when the areas $2\pi\md{a}$ of the holomorphic discs
in Lemma \ref{fi5lem} with boundary in the fibres of $F,F'$
shrink to zero, at $a=0$. This should also be something which
happens in SL fibrations of almost Calabi--Yau 3-folds, when
holomorphic discs with boundaries in the fibres shrink to a point.

\section{A model for codimension 2 singular behaviour}
\label{fi6}

Theorems \ref{fi5thm1} and \ref{fi5thm2} modelled the singular
behaviour the author expects to occur in codimension one in
generic SL fibrations. We shall now construct a model for
the next most generic kind of singular behaviour, which
occurs in codimension two in generic SL fibrations.

Unfortunately we cannot write the fibration down explicitly,
so we will construct it using the analytic results of \S\ref{fi4}
and describe its properties. We begin by defining a family of
solutions $(\hat u_{a,\al},\hat v_{a,\al})$ of~\eq{fi4eq5}.

\begin{definition} Let $D$ be the unit disc $\bigl\{(x,y)\in\R^2:x^2+
y^2\le 1\bigr\}$ in $\R^2$, with boundary ${\cal S}^1$, the
unit circle. Define a coordinate $\th:\R/2\pi\Z\ra{\cal S}^1$
by $\th\mapsto(\cos\th,\sin\th)$. Then $\cos(j\th),\sin(j\th)\in
C^\iy({\cal S}^1)$ for~$j\ge 1$.

For all $a,\al\in\R$ let $\hat f_{a,\al}$ be the unique solution
of \eq{fi4eq12} in $D$ (with weak derivatives) given by Theorem
\ref{fi4thm2}, with this value of $a$ and boundary condition
\e
\hat f_{a,\al}\vert_{{\cal S}^1}=\hat\phi_\al=\al\cos\th-\cos(3\th).
\label{fi6eq1}
\e
Then $\hat f_{a,\al}\in C^\iy(D)$ if $a\ne 0$ and $\hat f_{0,\al}
\in C^1(D)$. Define $\hat u_{a,\al}=\frac{\pd}{\pd y}\hat f_{a,\al}$
and $\hat v_{a,\al}=\frac{\pd}{\pd x}\hat f_{a,\al}$. Then
$(\hat u_{a,\al},\hat v_{a,\al})$ is a solution of \eq{fi4eq5}
in $C^\iy(D)$ if $a\ne 0$, and a singular solution of \eq{fi4eq4}
in $C^0(D)$ if $a=0$, in the sense of Definition~\ref{fi4def1}.
\label{fi6def1}
\end{definition}

\begin{theorem} These solutions $(\hat u_{a,\al},\hat v_{a,\al})$
have the following properties:
\begin{itemize}
\setlength{\itemsep}{0pt}
\setlength{\parsep}{0pt}
\item[{\rm(i)}] $\hat u_{a,\al}(x,-y)\equiv -\hat u_{a,\al}(x,y)$
and\/ $\hat v_{a,\al}(x,-y)\equiv \hat v_{a,\al}(x,y)$.
\item[{\rm(ii)}] $\hat u_{a,\al}(-x,y)\equiv -\hat u_{a,\al}(x,y)$
and\/ $\hat v_{a,\al}(-x,y)\equiv\hat v_{a,\al}(x,y)$.
\item[{\rm(iii)}] $\hat u_{-a,\al}(x,y)\equiv\hat u_{a,\al}(x,y)$
and\/ $\hat v_{-a,\al}(x,y)\equiv\hat v_{a,\al}(x,y)$.
\item[{\rm(iv)}] There exists $C>0$ with\/ $\md{\hat u_{a,\al}}\le C$
and\/ $\md{\hat v_{a,\al}-\al}\le C$ on $D$ for all\/~$a,\al$.
\item[{\rm(v)}] If\/ $a\in\R$, $(b,c)\in D^\circ$ and\/
$\al<\al'$ then $\hat v_{a,\al}(b,c)<\hat v_{a,\al'}(b,c)$.
\item[{\rm(vi)}] For all $a,\al\in\R$ the function
$x\mapsto\hat v_{a,\al}(x,0)$ is strictly increasing for
$x\in[-1,0]$, and strictly decreasing for~$x\in[0,1]$.
\item[{\rm(vii)}] For all $a,\al\in\R$ the function
$y\mapsto\hat v_{a,\al}(0,y)$ is strictly decreasing for
$y\in[-1,0]$, and strictly increasing for~$y\in[0,1]$.
\end{itemize}
\label{fi6thm1}
\end{theorem}

\begin{proof} As $\hat\phi_\al$ has the symmetries
$\hat\phi_\al(x,-y)=-\hat\phi_\al(-x,y)=\hat\phi_\al(x,y)$,
uniqueness in Theorem \ref{fi4thm2} gives $\hat f_{a,\al}(x,-y)
=-\hat f_{a,\al}(-x,y)=\hat f_{a,\al}(x,y)$, and parts (i), (ii)
follow by taking partial derivatives. Part (iii) holds as
\eq{fi4eq5} depends only on $a^2$ rather than $a$, and
$\hat\phi_\al$ is independent of~$a$.

An important part of the proof of Theorem \ref{fi4thm2} in
\cite{Joyc2,Joyc3} was to derive {\it a priori estimates}
for $\cnm{u}{0}$ and $\cnm{v}{0}$ in terms of $\cnm{\phi}{2}$.
This was done by using functions of the form $\be x+\ga y+\de$
as super- and subsolutions for $f$ at each point of $\pd S$,
and so derive a bound for $\md{\pd f}$ on $\pd S$. As the
maxima of $u,v$ occur on $\pd S$, this implies bounds for
$\md{u},\md{v}$ on~$S$.

Now in our case this can be done uniformly in $\al$. That
is, if $\be x+\ga y+\de\le -\cos(3\th)$ on ${\cal S}^1$
then $(\al+\be)x+\ga y+\de\le\hat\phi_\al$ on ${\cal S}^1$
for all $\al$, and therefore $(\al+\be)x+\ga y+\de\le\hat
f_{a,\al}$ on $D$ for all $a,\al$, and a similar statement for
supersolutions. Following the proof of \cite[Th.~3.9]{Joyc2},
we easily deduce part~(iv).

Suppose $a\in\R$, $(b,c)\in D^\circ$, $\al<\al'$ and
$\hat v_{a,\al}(b,c)=\hat v_{a,\al'}(b,c)$. Define
\begin{align*}
u_1(x,y)&=\hat u_{a,\al}(x,y)+\hat u_{a,\al'}(b,c)-\hat u_{a,\al}(b,c),
\;\>& v_1(x,y)&=\hat v_{a,\al}(x,y),\\
u_2(x,y)&=\hat u_{a,\al'}(x,y)
\qquad\qquad\qquad\text{and}\quad&
v_2(x,y)&=\hat v_{a,\al'}(x,y).
\end{align*}
Then $(u_j,v_j)$ satisfy \eq{fi4eq5} in $D$, and $(u_1,v_1)=(u_2,v_2)$
at~$(b,c)\in D^\circ$.

Furthermore, $(u_1,v_1)$, $(u_2,v_2)$ come from Theorem
\ref{fi4thm2} with boundary data
\begin{equation*}
\phi_1\!=\!\al\cos\th\!-\!\cos(3\th)\!+\!
\bigl(\hat u_{a,\al'}(b,c)\!-\!\hat u_{a,\al}(b,c)\bigr)\sin\th
\;\>\text{and}\;\>
\phi_2\!=\!\al'\cos\th\!-\!\cos(3\th).
\end{equation*}
Thus $\phi_1-\phi_2$ is a nontrivial linear combination of
$\cos\th,\sin\th$, and has exactly 1 local maximum and 1 local
minimum on ${\cal S}^1$. Applying Theorem \ref{fi4thm4} with
$l=1$ shows that $(u_1,v_1)-(u_2,v_2)$ has no zeroes in $D^\circ$.
But this contradicts $(u_1,v_1)=(u_2,v_2)$ at~$(b,c)$.

So given $a\in\R$ and $(b,c)\in D^\circ$, whenever $\al<\al'$
we have $\hat v_{a,\al}(b,c)\ne\hat v_{a,\al'}(b,c)$. Since
$\hat v_{a,\al}(b,c)$ depends continuously on $\al,a,b,c$ it follows
that $v_{a,\al}(b,c)$ is either a strictly increasing or a strictly
decreasing function of $\al$, and which of the two is independent
of $a,b,c$. Then (iv) shows it is increasing, proving~(v).

For $\be\in\R$, we can consider $(0,\be)$ to be a solution
of \eq{fi4eq5} on $D$, coming from Theorem \ref{fi4thm2}
with $f=\be x$ and boundary data $\phi=\be\cos\th$. Now
$\al\cos\th-\cos(3\th)-\be\cos\th$ has at most 3 local maxima
and 3 local minima on ${\cal S}^1$ by \cite[Prop.~10.2]{Joyc4}.
Hence, applying Theorem \ref{fi4thm4} with $l=3$ shows that
$(\hat u_{a,\al},\hat v_{a,\al})-(0,\be)$ has at most two zeroes
in $D^\circ$, for any~$\be\in\R$.

We shall use this to show that
\begin{itemize}
\setlength{\itemsep}{0pt}
\setlength{\parsep}{0pt}
\item[(a)] Suppose $\md{x},\md{x'}<1$ and
$\hat v_{a,\al}(x,0)=\hat v_{a,\al}(x',0)$. Then~$x=\pm x'$.
\item[(b)] Suppose $0<\md{x},\md{y}<1$.
Then~$\hat v_{a,\al}(x,0)\ne\hat v_{a,\al}(0,y)$.
\end{itemize}
Part (i) gives $\hat u_{a,\al}(x,0)\equiv 0$. So if $\md{x},\md{x'}<1$
and $\hat v_{a,\al}(x,0)=\hat v_{a,\al}(x',0)=\be$ then
$(\hat u_{a,\al},\hat v_{a,\al})=(0,\be)$ at $(\pm x,0)$ and
$(\pm x',0)$, by part (ii). As there are at most two such points
we must have $x=\pm x'$, proving~(a).

Similarly, (ii) gives $\hat u_{a,\al}(0,y)\equiv 0$. Thus if
$0<\md{x},\md{y}<1$ and $\hat v_{a,\al}(x,0)=\hat v_{a,\al}(0,y)=\be$
then by (i), (ii) we see that $(\hat u_{a,\al},\hat v_{a,\al})=(0,\be)$
at the four points $(\pm x,0)$ and $(0,\pm y)$, a contradiction.
This proves~(b).

Now from (a), (b) and the continuity of $\hat v_{a,\al}$ it is not
difficult to see that {\it either} (vi), (vii) hold, {\it or}
(vi), (vii) hold, but swapping `increasing' and `decreasing'
throughout. But
\begin{equation*}
\int_{-1}^1\hat v_{a,\al}(x,0)\,\d x
=\hat\phi_\al(1,0)-\hat\phi_\al(-1,0)=2\al-2,
\end{equation*}
so the average of $\hat v_{a,\al}$ on the $x$-axis is $\al-1$, and
\begin{equation*}
\hat v_{a,\al}(0,1)=\frac{\pd\hat f_{a,\al}}{\pd x}(0,1)=
-\frac{\d\hat\phi_\al}{\d\th}(0,1)=\al+3.
\end{equation*}
Therefore $\hat v_{a,\al}$ is greater on the $y$-axis than on the
$x$-axis, and so (vi) and (vii) hold, rather than their opposites.
\end{proof}

Next we identify the singularities of~$(\hat u_{0,\al},\hat v_{0,\al})$.

\begin{proposition} There exist unique $\al_0<\al_1$ in $\R$ such that:
\begin{itemize}
\setlength{\itemsep}{0pt}
\setlength{\parsep}{0pt}
\item[{\rm(i)}] If\/ $\al\notin[\al_0,\al_1]$ then
$(\hat u_{0,\al},\hat v_{0,\al})$ has no singularities in~$D$.
\item[{\rm(ii)}] $(\hat u_{0,\al_0},\hat v_{0,\al_0})$ has a
singularity of multiplicity $2$ and maximum type at\/ $(0,0)$,
and no other singularities.
\item[{\rm(iii)}] If\/ $\al\in(\al_0,\al_1)$ there exists
$x\in(0,1)$ such that\/ $(\hat u_{0,\al_0},\hat v_{0,\al_0})$ has
a singularity of multiplicity $1$ and increasing type at\/ $(-x,0)$,
a singularity of multiplicity $1$ and decreasing type at\/ $(x,0)$,
and no other singularities.
\item[{\rm(iv)}] $(\hat u_{0,\al_1},\hat v_{0,\al_1})$ is
singular at\/ $(\pm 1,0)$ on $\pd D$, and has no other
singularities.
\end{itemize}
\label{fi6prop}
\end{proposition}

\begin{proof} This follows quickly from Theorem \ref{fi6thm1} using
the Intermediate Value Theorem, except for the multiplicities of
the singular points. To find these we apply Theorem \ref{fi4thm6} to
$(\hat u_{0,\al},\hat v_{0,\al})$. We have $\phi-\phi'=2\al\cos\th
-2\cos(3\th)$, with at most 3 local maxima and 3 local minima
on ${\cal S}^1$ by \cite[Prop.~10.2]{Joyc4}. So Theorem \ref{fi4thm6}
shows that there are at most two singularities of $(\hat u_{0,\al},
\hat v_{0,\al})$ in $D^\circ$, counted with multiplicity.
Multiplicity 1 in (iii) follows at once, and multiplicity 2 in (ii)
from Proposition \ref{fi4prop5}, as $(0,0)$ is of maximum type.
\end{proof}

Next we use the results of \S\ref{fi45} to construct the
special Lagrangian fibration we want. We apply Example
\ref{fi4ex2} with $S=D$ and $\phi=-\cos(3\th)$. Equivalently,
we apply Definition \ref{fi4def2} and Theorem \ref{fi4thm3}
with $S=D$ and $\Phi(a,\al,\be)=\al\cos\th+\be\sin\th-\cos(3\th)$.
Writing the definition and theorem out explicitly in our case gives:

\begin{definition} For each $a,\al,\be\in\R$, define $\hat N_{a,\al,\be}$
in $\C^3$ by
\e
\begin{split}
\hat N_{a,\al,\be}=\bigl\{(z_1&,z_2,z_3)\in\C^3:
\ms{z_1}-\ms{z_2}=2a,\;\> x,y\in\R,\;\> x^2+y^2<1,\\
&z_1z_2=\hat v_{a,\al}(x,y)+iy,\;\>
z_3=x+i\hat u_{a,\al}(x,y)+i\be\bigr\}.
\end{split}
\label{fi6eq2}
\e
Then $\hat N_{a,\al,\be}$ is a noncompact SL 3-fold without
boundary in~$\C^3$.
\label{fi6def2}
\end{definition}

\begin{theorem} In the situation above, distinct\/ $\hat N_{a,\al,\be}$
are disjoint. Define
\e
V=\bigl\{(z_1,z_2,z_3)\in\C^3:(\Re z_3)^2+(\Im z_1z_2)^2<1\bigr\}.
\label{fi6eq3}
\e
Then there exists a continuous, surjective $\hat F:V\ra\R^3$
with\/ $\hat F^{-1}(a,\al,\be)=\hat N_{a,\al,\be}$ for all\/
$(a,\al,\be)\in\R^3$. Thus, $\hat F$ is a special Lagrangian
fibration of\/~$V$.
\label{fi6thm2}
\end{theorem}

Proposition \ref{fi6prop} gives the {\it discriminant} of $\hat F$.

\begin{corollary} In the situation above, the discriminant of\/ $\hat F$ is
\e
\hat\De=\bigl\{(0,\al,\be):\al\in [\al_0,\al_1),\quad\be\in\R\bigr\},
\label{fi6eq4}
\e
and the set of singular points is $\bigl\{(0,0,z_3):z_3\in\C$,
$\md{\Re z_3}<1\bigr\}$.
\label{fi6cor}
\end{corollary}

\subsection{Discussion}
\label{fi61}

We have constructed an SL fibration $\hat F$ in which the
discriminant $\hat\De$ is a {\it ribbon}, a portion of a plane in
$\R^3$. As in \S\ref{fi5}, $\hat\De$ is of {\it real codimension
one} in the base $\R^3$. However, we are particularly interested in
what happens at the boundary $\bigl\{(0,\al_0,\be):\be\in\R\bigr\}$
of $\hat\De$, which is of {\it real codimension two} in~$\R^3$.

In the interior of $\hat\De$, each singular fibre $\hat N_{0,\al,\be}$
with $\al\in(\al_0,\al_1)$ has two singular points, both of multiplicity
1 and one each of increasing and decreasing type. Locally these
singularities are modelled on SL $T^2$-cones like $C$ in \eq{fi5eq12}.
As $\al\ra\al_0$, these two singular points come together, until at
$\al=\al_0$ they fuse to form a different kind of singularity, of
multiplicity 2. For $\al<\al_0$ there are no singularities. Thus
the picture is that as $\al$ decreases through $\al_0$, the two
singular points in $\hat N_{0,\al,\be}$ come together and cancel
out. As $\al$ increases through $\al_1$ the singular points cross
the boundary of~$V$.

Now consider the set of singular points $\bigl\{(0,0,z_3):z_3\in\C$,
$\md{\Re z_3}<1\bigr\}$. When $\Re z_3\in(-1,0)$, we see from part
(iii) of Proposition \ref{fi6prop} that the fibre of $\hat F$ passing
through $(0,0,z_3)$ has a singularity of multiplicity 1 and increasing
type, like those of $F$ in Theorem \ref{fi5thm1}. When $\Re z_3\in
(0,1)$, the fibre of $\hat F$ through $(0,0,z_3)$ has a singularity
of multiplicity 1 and decreasing type, like those of $F'$ in
Theorem~\ref{fi5thm2}.

So the situation is that near $\bigl\{(0,0,z_3):z_3\in\C$,
$\Re z_3\in(-1,0)\bigr\}$, the fibration $\hat F$ is locally
modelled on $F$ in Theorem \ref{fi5thm1}, and near $\bigl\{
(0,0,z_3):z_3\in\C$, $\Re z_3\in(0,1)\bigr\}$, the fibration
$\hat F$ is locally modelled on $F'$ in Theorem \ref{fi5thm2}.
On the line $\bigl\{(0,0,z_3):z_3\in\C$, $\Re z_3=0\bigr\}$,
we have multiplicity 2 singularities, which mark the transition
between the $F$ and $F'$ local models.

I claim that $\hat F$ is a local model for generic SL fibrations.

\begin{conjecture} The SL fibration $\hat F$ of Theorem \ref{fi6thm2}
is a generic local model for one kind of codimension two
singularities of SL fibrations $f:M\ra B$ of almost Calabi--Yau
$3$-folds, using `generic' in the sense of Definition~\ref{fi3def4}.
\label{fi6conj}
\end{conjecture}

Again, I am not defining what I mean by a `local model' here,
but basically the fibrations should have the same topological
structure and the same kinds of singularities. From above, for
$\al\in(\al_0,\al_1)$ the singular fibres $\hat N_{0,\al,\be}$
have two singular points modelled on the SL $T^2$-cone $C$ in
\eq{fi5eq12}. Thus it is an SL 3-fold with {\it conical
singularities}, as in Definition~\ref{fi2def6}.

Therefore, as in \S\ref{fi51}, we can use the material of
\S\ref{fi23} to give a partial proof of the genericity of
$\hat F$ around the fibres $\hat N_{0,\al,\be}$ for
$\al\in(\al_0,\al_1)$. However, the singularities of
$\hat N_{0,\al_0,\be}$ are {\it not\/} conical in the
sense of Definition \ref{fi2def6}, so this approach will
not help around~$\al=\al_0$.

We shall discuss one feature of these fibrations, and why it is
generic, in more detail. For $\al\in(\al_0,\al_1)$ the fibre
$\hat N_{0,\al,\be}$ has {\it two} singular points. {\it A priori}
this seems unlikely: distinct singular points ought to occur
independently, and so for a fibre to have two codimension one
singular points should be a codimension two phenomenon, not
codimension one as in $\hat F$. However, this is not the case.

We can explain this in terms of {\it holomorphic discs}, as in
\S\ref{fi51}. Let us identify the holomorphic discs in $\C^3$
with boundary in~$\hat N_{a,\al,\be}$.

\begin{lemma} Let\/ $a>0$, $\al,\be\in\R$ and\/ $x\in(0,1)$ with\/
$\hat v_{a,\al}(x,0)=0$. Then
\begin{itemize}
\setlength{\itemsep}{0pt}
\setlength{\parsep}{0pt}
\item[{\rm(a)}] $D_\pm=\bigl\{(z_1,0,\pm x+i\hat u_{a,\al}(x,0)
+i\be):z_1\in\C$, $\md{z_1}^2\le 2a\bigr\}$ are two holomorphic
discs with boundary in $\hat N_{a,\al,\be}$ and area $2\pi a$, and
\item[{\rm(b)}] $D'_\pm=\bigl\{(0,z_2,\pm x+i\hat u_{a,\al}(x,0)
+i\be):z_2\in\C$, $\md{z_2}^2\le 2a\bigr\}$ are two holomorphic
discs with boundary in $\hat N_{-a,\al,\be}$ and area~$2\pi a$.
\end{itemize}
\label{fi6lem}
\end{lemma}

The proof is trivial. Suppose now that $a>0$ is small,
$\al\in(\al_0,\al_1)$ and $\be\in\R$. Then by Theorem \ref{fi6thm1}
and Proposition \ref{fi6prop} there exists $x\in(0,1)$ with
$\hat v_{a,\al}(x,0)=0$, and Lemma \ref{fi6lem} gives holomorphic
discs $D_\pm$ with boundary in $\hat N_{a,\al,\be}$, and area
$2\pi a$. These discs are {\it homologous} in $H_2(\C^3,\hat
N_{a,\al,\be};\R)$. As $a$ decreases to zero $D_\pm$ collapse to
points, and $\hat N_{a,\al,\be}$ develops two singular points.

Recall from \S\ref{fi51} that if $D$ is a holomorphic disc in an
almost Calabi--Yau 3-fold $(M,J,\om,\Om)$ with boundary in an SL
3-fold $L$ in $M$, then the area of $D$ is $[\om]\cdot[D]$, where
$[D]\in H_2(M,L;\Z)$ is the relative homology class of $D$.
Therefore, if holomorphic discs $D_1,\ldots,D_k$ with boundary
in $L$ have the same homology class in $H_2(M,L;\Z)$ then they
have the same area. If we deform $L$ so that this area becomes
zero, then $D_1,\ldots,D_k$ will simultaneously collapse to
points, and $L$ will develop $k$ singular points. For some
rigorous results on this when $k=2$, see~\cite[\S 10.3]{Joyc10}.

This shows that distinct singular points of SL 3-folds may not
be independent. Instead, if a singular point results from the
collapse of a holomorphic disc, then singular points from the
collapse of homologous holomorphic discs will always occur
together. This is why it is permissible for the fibres of
$\hat F$ to have two singular points in codimension one,
and for this to be generic.

We can also use this to explain why the discriminant $\hat\De$
can have a boundary. For $\al\in(\al_0,\al_1)$ and $a$ small,
there are two holomorphic discs $D_\pm$ with boundary in
$N_{a,\al,\be}$. These discs have {\it opposite sign}, so
that the number of holomorphic discs counted with signs is
zero. As $\al$ decreases with $a$ fixed, it reaches a value
$\al'\approx\al_0$ with $v_{a,\al'}(0,0)=0$, and then $D_\pm$
come together and cancel out. For $\al<\al'$ there are no
holomorphic discs with boundary in~$N_{a,\al,\be}$.

Thus, if we pass through the hypersurface $a=0$ when
$\al\in(\al_0,\al_1)$, two holomorphic discs collapse
to two singularities. But if we decrease $a$ past $\al_0$,
the two holomorphic discs cancel, and then we can pass
through $a=0$ without a singularity, as there are no
holomorphic discs to collapse. 

\section{\!\!How smooth SL fibrations become non-smooth}
\label{fi7}

Next we shall extend the results of \S\ref{fi4} from strictly
convex domains in $\R^2$ to solutions on a strip in $\R^2$
which are periodic under a group of translations. We shall
use this to model what happens near an ${\cal S}^1$ singularity
of a singular fibre of a smooth SL fibration, and so describe
how smooth SL fibrations become non-smooth under small deformations.

\subsection{A class of periodic $\U(1)$-invariant SL 3-folds}
\label{fi71}

Here is the situation we shall work with.

\begin{definition} Let $R,P>0$, and define $S=\bigl\{(x,y)\in\R^2:
\md{y}\le R\bigr\}$. We shall study functions $u,v:S\ra\R$
which satisfy \eq{fi4eq4} or \eq{fi4eq5}, and the
periodicity condition
\e
u(x+P,y)\equiv u(x,y) \quad\text{and}\quad
v(x+P,y)\equiv v(x,y) \quad\text{for all $(x,y)\in S$.}
\label{fi7eq1}
\e
\label{fi7def1}
\end{definition}

The main point is that although $S$ is noncompact, equation
\eq{fi7eq1} implies that $u,v$ are invariant under the
$\Z$-action
\e
(x,y)\,{\buildrel n\over\longmapsto}\,(x+nP,y)
\quad\text{for $n\in\Z$,}
\label{fi7eq2}
\e
so we can treat $u,v$ as functions on the annulus $S/\Z$.
Since $S/\Z$ is compact, analytic methods which rely on
compactness will still apply, such as existence results
for the Dirichlet problem.

Here are two lemmas on the consequences of the periodicity
conditions \eq{fi7eq1}. The second is an analogue of
Proposition~\ref{fi4prop2}.

\begin{lemma} Let\/ $R,P$ and\/ $S$ be as above, $a\ne 0$, and\/
$u,v\in C^1(S)$ satisfy \eq{fi4eq5} and\/ \eq{fi7eq1}. Then
there is a unique $\ga\in\R$ with
\e
\int_0^Pv(x,y)\,\d x=\ga P\quad\text{for all $y\in[-R,R]$.}
\label{fi7eq3}
\e
\label{fi7lem1}
\end{lemma}

\begin{proof} Since $u,v$ are continuously differentiable and
satisfy \eq{fi4eq5}, we have
\begin{equation*}
\frac{\d}{\d y}\int_0^Pv(x,y)\d x=\!
\int_0^P\frac{\pd v}{\pd y}(x,y)\d x=\!
\int_0^P\frac{\pd u}{\pd x}(x,y)\d x=u(P,y)-u(0,y)=0,
\end{equation*}
using \eq{fi7eq1}. Hence $\int_0^Pv(x,y)\,\d x$ is independent
of $y$, and \eq{fi7eq3} holds for some unique~$\ga$.
\end{proof}

\begin{lemma} Let\/ $R,P$ and\/ $S$ be as above, and\/ $v\in C^2(S)$
satisfy \eq{fi4eq6} for $a\ne 0$, $v(x+P,y)\equiv v(x,y)$ for all\/
$(x,y)\in S$, and\/ $\int_0^Pv(x,R)\,\d x=\int_0^Pv(x,-R)\,\d x$.
Then there exists a unique $u\in C^2(S)$ with\/ $u(0,0)=0$ such
that\/ $u,v$ satisfy \eq{fi4eq5} and\/ $u(x+P,y)\equiv u(x,y)$
for all\/~$(x,y)\in S$.
\label{fi7lem2}
\end{lemma}

\noindent{\it Proof.} The proof of Proposition \ref{fi4prop2} shows
that there exists $u\in C^2(S)$, unique up to $u\mapsto u+c$, such
that $u,v$ satisfy \eq{fi4eq5}. Requiring $u(0,0)=0$ fixes $c$. As
$v$ is periodic under \eq{fi7eq2} we see from \eq{fi4eq5} that $\pd
u$ is periodic, and hence $u$ satisfies $u(x+P,y)\equiv u(x,y)+\de$
for some $\de\in\R$. We must show $\de=0$. As $\frac{\pd u}{\pd x}
=\frac{\pd v}{\pd y}$, this follows from
\begin{align*}
0&=\int_0^Pv(x,R)\,\d x-\int_0^Pv(x,-R)\,\d x=
\int_0^P\int_{-R}^R\frac{\pd v}{\pd y}(x,y)\,\d y\,\d x\\
&=\int_{-R}^R\int_0^P\frac{\pd u}{\pd x}(x,y)\,\d x\,\d y
=\int_{-R}^R\bigl(u(P,y)-u(0,y)\bigr)\d y=2R\de.
\tag*{$\square$}
\end{align*}

In \S\ref{fi6} our main tool for constructing SL fibrations
was Theorem \ref{fi4thm2}, the Dirichlet problem for solutions
$f$ of \eq{fi4eq12}. We could take the same approach in this
situation, but it turns out to be more elegant to solve the
Dirichlet problem for solutions $v$ of \eq{fi4eq6} instead.

Now \cite[Th.~8.8]{Joyc2} proves existence and uniqueness
for the Dirichlet problem for $v$ in \eq{fi4eq6} on domains
in $\R^2$ when $a\ne 0$. The proof does not assume the domain
is convex, or use the fact that domains are contractible. In
fact the proof applies to the compact annulus $S/\Z$ without
change. This gives:

\begin{theorem} Let\/ $R,P$ and\/ $S$ be as above, $a\ne 0$, $k\ge 0$,
and\/ $\al\in(0,1)$. Suppose $\phi^\pm\in C^{k+2,\al}(\R)$ with\/
$\phi^\pm(x+P)\equiv\phi^\pm(x)$. Then there exists a unique
$v\in C^{k+2,\al}(S)$ satisfying \eq{fi4eq6},
$v(x+P,y)\equiv v(x,y)$ and\/~$v(x,\pm R)\equiv\phi^\pm(x)$.
\label{fi7thm1}
\end{theorem}

Before we extend this to the case $a=0$, we prove some a priori
estimates for solutions $u,v$ in terms of the boundary data
$\phi^\pm$ for~$v$.

\begin{theorem} Let\/ $R,P$ and\/ $S$ be as above, $a\ne 0$ and\/
$u,v\in C^2(S)$ satisfy \eq{fi4eq5}, \eq{fi7eq1}, and\/ $u(0,0)=0$.
Define $\phi^\pm\in C^2(\R)$ by $v(x,\pm R)=\phi^\pm(x)$. Then
\e
\cnm{v}{0}\!=\!\max\bigl(\cnm{\phi^+}{0},\cnm{\phi^-}{0}\bigr)
\;\>\text{and}\;\>
\bcnm{\ts\frac{\pd v}{\pd x}}{0}\!=\!\max\bigl(
\bcnm{\ts\frac{\pd\phi^+}{\pd x}}{0},
\bcnm{\ts\frac{\pd\phi^-}{\pd x}}{0}\bigr).
\label{fi7eq4}
\e
Suppose further that\/ $A,B>0$ with\/ $\md{a}\le A$, and\/
$\cnm{\phi^\pm}{2}\le B$. Then there exist $C,D>0$ depending
only on $R,P,A,B$ such that\/ $\cnm{u}{0}\le C$,
$\cnm{\pd u\vert_{\pd S}}{0}\le D$
and\/~$\cnm{\pd v\vert_{\pd S}}{0}\le D$.
\label{fi7thm2}
\end{theorem}

\begin{proof} Regard $v$ and $\frac{\pd v}{\pd x}$ as functions
on the compact annulus $S/\Z$, as they are invariant under the
$\Z$-action \eq{fi7eq2}. Now using the {\it maximum principle}
for elliptic equations of a certain form, \cite[Cor.~4.4]{Joyc2}
shows that the maximum and minimum of $v$ on a domain $T$ occur
on $\pd T$, and \cite[Prop.~8.12]{Joyc2} that the maximum of
$\bmd{\frac{\pd v}{\pd x}}$ occurs on $\pd T$. These proofs are
also valid on the compact annulus $S/\Z$. But $v(x,\pm R)=
\phi^\pm(x)$ and $\frac{\pd v}{\pd x}(x,\pm R)=\frac{\pd}{\pd x}
\phi^\pm(x)$ on $\pd(S/\Z)$. Maximizing and minimizing then
gives~\eq{fi7eq4}.

Next we estimate $\md{\pd u},\md{\pd v}$ on $\pd S$, following
\cite[Prop.~8.6]{Joyc3}. Observe that if $(x,y)\in S$ with
$\md{y}\ge\ha R$ then
\begin{equation*}
\ts\frac{1}{4}R^2\le v(x,y)^2+y^2+a^2\le B^2+R^2+A^2,
\end{equation*}
since $\md{v}\le\max\bigl(\cnm{\phi^+}{0},\cnm{\phi^-}{0}\bigr)\le B$
by \eq{fi7eq4} and $\cnm{\phi^\pm}{2}\le B$. It follows that we can
treat \eq{fi4eq6} as a quasilinear elliptic equation on $v$, which
is {\it uniformly elliptic} on $\md{y}\ge\ha R$, with constants
of ellipticity depending only on $A,B$ and~$R$.

Now Gilbarg and Trudinger \cite[Th.~14.1, p.~337]{GiTr} show that
if $v\in C^2(T)$ satisfies a quasilinear equation $Qv=0$ of the
form \eq{fi4eq2} on a domain $T$ and $v\vert_{\pd T}=\phi\in
C^2(\pd T)$, then $\cnm{\pd v\vert_{\pd T}}{0}\le K$ for some
$K>0$ depending only on $T$, upper bounds for $\cnm{v}{0}$ and
$\cnm{\phi}{2}$, and certain constants to do with $Q$, which
ensure that $Q$ is uniformly elliptic and $b$ not too large.

As this is a local result, it is enough for the conditions to
hold within distance $\ha R$ of $\pd T$. So there exists $K>0$
depending only on $A,B$ and $R$ such that $\cnm{\pd v
\vert_{\pd S}}{0}\le K$. But $\pd v$ determines $\pd u$ by
\eq{fi4eq5}, and we easily deduce a bound for $\cnm{\pd u
\vert_{\pd S}}{0}$. Thus there exists $D>0$ depending only
on $A,B$ and $R$ such that $\cnm{\pd u\vert_{\pd S}}{0}\le D$
and~$\cnm{\pd v\vert_{\pd S}}{0}\le D$.

Finally we estimate $\cnm{u}{0}$. Since $u(x,R)$ is periodic in
$x$ with period $P$ and $\bmd{\frac{\pd u}{\pd x}(x,R)}\le D$,
we see that the {\it variation} of $x\mapsto u(x,R)$ is at most
$\ha PD$. Similarly, the variation of $x\mapsto u(x,-R)$ is at
most $\ha PD$. Therefore $u$ varies only a bounded amount from
its average value on $y=R$, and the same on $y=-R$. We need to
bound the difference between these average values. To do this
we use the method of \cite[\S 8.4]{Joyc3} to bound
$\lnm{\frac{\pd u}{\pd y}}{1}$ on~$S/\Z$.

Define $J(a,v)=-\int_0^v(w^2+a^2)^{-1/4}\d w$. Then
\cite[Prop.~8.9]{Joyc3} shows that
\e
\begin{split}
&\int_0^P\int_{-R}^R(v^2+a^2)^{-1/4}
\Bigl[\ha(v^2+y^2+a^2)^{-1/2}\Bigl(\frac{\pd v}{\pd x}\Bigr)^2
+\Bigl(\frac{\pd v}{\pd y}\Bigr)^2\Bigr]\d y\,\d x=\\
&\int_0^PJ(a,v(x,R))\frac{\pd u}{\pd x}(x,R)\,\d x
-\int_0^PJ(a,v(x,-R))\frac{\pd u}{\pd x}(x,-R)\,\d x.
\end{split}
\label{fi7eq5}
\e
The proof is that $\int_{S/\Z}\d\bigl(J(a,v)\d u\bigr)=
\int_{\pd(S/\Z)}J(a,v)\d u$ by Stokes' Theorem, and using
\eq{fi4eq5} to rewrite the l.h.s.\ gives \eq{fi7eq5}. Therefore
\e
\begin{split}
&\Bigl(\int_0^P\int_{-R}^R
\Big\vert\frac{\pd u}{\pd y}\Big\vert\d y\,\d x\Bigr)^2
\Bigl(\int_0^P\int_{-R}^R(v^2+y^2+a^2)^{-1/4}\d y\,\d x\Bigr)^{-1}\\
\le\,&\int_0^P\int_{-R}^R
(v^2+y^2+a^2)^{1/4}\Bigl(\frac{\pd u}{\pd y}\Bigr)^2\d y\,\d x\\
\le\,&\ha\int_0^P\!\!J\bigl(a,v(x,R)\bigr)\frac{\pd u}{\pd x}(x,R)\d x
\!-\!\ha\int_0^P\!\!J\bigl(a,v(x,-R)\bigr)\frac{\pd u}{\pd x}
(x,-R)\d x,
\end{split}
\label{fi7eq6}
\e
where the second line follows from H\"older's inequality,
and the third from \eq{fi7eq5} and the inequality
$(v^2+y^2+a^2)^{1/4}\ms{\frac{\pd u}{\pd y}}\le
\frac{1}{4}(v^2+a^2)^{-1/4}(v^2+y^2+a^2)^{-1/2}
\ms{\frac{\pd v}{\pd x}}$, which follows from~\eq{fi4eq5}.

Using $\md{a}\le A$, $\md{v}\le B$ and the methods of
\cite[Prop.~8.10]{Joyc3} we may derive a priori estimates
depending only on $A,B,R$ and $P$ for the third line of
\eq{fi7eq6} and the second integral on the first line.
This gives an upper bound, $E$ say, for
$\int_0^P\int_{-R}^R\bmd{\frac{\pd u}{\pd y}}\,\d y\,\d x$. But
\begin{equation*}
\Big\vert\int_0^P\bigl(u(x,R)\!-\!u(x,-R)\bigr)\,\d x\Big\vert\!=\!
\Big\vert\int_0^P\!\!\!\!\int_{-R}^R\frac{\pd u}{\pd y}\,\d y\,\d x
\Big\vert\le\int_0^P\!\!\!\!\int_{-R}^R\Big\vert\frac{\pd u}{\pd y}
\Big\vert\,\d y\,\d x\!\le\!E.
\end{equation*}
Therefore the difference in the average values of $u$ on
$y=\pm R$ is at most $E/P$. But the variation of $u$ on $y=R$
and on $y=-R$ is at most $\ha PD$, from above. It follows
that the variation of $u$ on both lines $y=\pm R$ together
is at most $C=E/P+PD$, which depends only on $A,B,R$ and~$P$.

Now $u$ satisfies a maximum principle by \cite[Cor.~4.4]{Joyc2},
and so the maximum and minimum of $u$ on the compact annulus
$S/\Z$ occur on $\pd(S/\Z)$. Hence the difference between the
maximum and minimum of $u$ is the variation of $u$ on both lines
$y=\pm R$ together, and is at most $C$. But $u(0,0)=0$, so the
maximum is nonnegative, and the minimum nonpositive. Therefore
the maximum is at most $C$ and the minimum at least $-C$, and
$\cnm{u}{0}\le C$, completing the proof.
\end{proof}

Our next result extends~\cite[Th.~8.17 \& 8.18]{Joyc3}.

\begin{theorem} Let\/ $R,P$ and\/ $S$ be as above, $a\in\R$,
$k\ge 0$, and\/ $\al\in(0,1)$. Suppose $\phi^+,\phi^-\in
C^{k+2,\al}(\R)$ with\/ $\phi^\pm(x+P)\equiv\phi^\pm(x)$
and\/ $\int_0^P\phi^+(x)\,\d x=\int_0^P\phi^-(x)\,\d x$.
Then if\/ $a\ne 0$ there exist unique $u,v\in C^{k+2,\al}(S)$
satisfying \eq{fi4eq5}, \eq{fi7eq1}, $u(0,0)=0$ and\/
$v(x,\pm R)\equiv \phi^\pm(x)$. If\/ $a=0$ there exist unique
$u,v\in C^0(S)$ which are a singular solution of\/ \eq{fi4eq4}
in the sense of Definition \ref{fi4def1}, and satisfy \eq{fi7eq1},
$u(0,0)=0$ and\/ $v(x,\pm R)\equiv \phi^\pm(x)$. Furthermore $u,v$
depend continuously in $C^0(S)$ on $\phi^+,\phi^-$ in
$C^{k+2,\al}(\R)$ and\/ $a$ in~$\R$.
\label{fi7thm3}
\end{theorem}

\begin{proof} When $a\ne 0$, existence and uniqueness of $v$
comes from Theorem \ref{fi7thm1}, and of $u$ from Lemma
\ref{fi7lem2}. To extend this to the case $a=0$ by taking the
limit $a\ra 0_+$ we follow \cite[\S 8]{Joyc3}, using the a
priori estimates of Theorem \ref{fi7thm2}. There are few
significant changes, and the problems caused by singular
points on the boundary in \cite{Joyc3} are absent in this
case, as there are no points $(x,0)$ on $\pd S$. The final
part is proved as in~\cite[Th.~8.18]{Joyc3}.
\end{proof}

Applying \cite[Prop.~8.7]{Joyc2} and \cite[Th.~6.16]{Joyc4}
on the annulus $S/\Z$ proves:

\begin{theorem} Let\/ $R,P$ and\/ $S$ be as above, $a\in\R$, $k\ge 0$
and\/ $\al\in(0,1)$, and suppose $\phi^\pm_i\in C^{k+2,\al}(\R)$
for $i=1,2$ satisfy $\phi^\pm_i(x+P)=\phi^\pm_i(x)$ and\/
$\phi^\pm_1(x)<\phi^\pm_2(x)$ for all\/ $x\in\R$, and\/
$\int_0^P\phi^+_i(x)\,\d x=\int_0^P\phi^-_i(x)\,\d x$. Let\/
$(u_i,v_i)$ be the (singular) solution of\/ \eq{fi4eq5}
produced in Theorem \ref{fi7thm3} from $\phi^\pm_i$ for
$i=1,2$. Then $v_1<v_2$ on~$S$.
\label{fi7thm4}
\end{theorem}

\subsection{Applications to SL fibrations}
\label{fi72}

We can now prove analogues of Definition \ref{fi4def2} and
Theorem~\ref{fi4thm3}.

\begin{definition} Let $R,P$ and $S$ be as above and $\al\in(0,1)$.
Suppose $\phi^+_{a,b},\phi^-_{a,b}\in C^{3,\al}(\R)$ are given
for all $a,b\in\R$ and satisfy
\begin{itemize}
\setlength{\itemsep}{0pt}
\setlength{\parsep}{0pt}
\item[(i)] $\phi^\pm_{a,b}$ depend continuously in $C^{3,\al}(\R)$
on~$a,b\in\R$,
\item[(ii)] $\phi^\pm_{a,b}(x+P)=\phi^\pm_{a,b}(x)$ for all~$a,b,x\in\R$,
\item[(iii)] $\int_0^P\phi^+_{a,b}(x)\,\d x=\int_0^P\phi^-_{a,b}(x)\,\d x$
for all~$a,b\in\R$,
\item[(iv)] If $b<b'$ then $\phi^\pm_{a,b}(x)<\phi^\pm_{a,b'}(x)$
for all~$a,x\in\R$.
\item[(v)] For all $a,x\in\R$ we have $\phi^+_{a,b}(x),
\phi^-_{a,b}(x)\ra\pm\iy$ as $b\ra\pm\iy$.
\end{itemize}
For all $a,b\in\R$, let $(u_{a,b},v_{a,b})$ be the (singular)
solution of \eq{fi4eq5} produced in Theorem \ref{fi7thm3} from
$\phi^\pm_{a,b}$. Let $\Z$ act on $\C^3$ by
\e
(z_1,z_2,z_3)\,{\buildrel n\over\longmapsto}\,(z_1,z_2,z_3+nP)
\quad\text{for $n\in\Z$,}
\label{fi7eq7}
\e
corresponding to \eq{fi7eq2}. For all $a,b,c\in\R$ define
$N_{a,b,c}$ in $\C^3/\Z$ by
\e
\begin{split}
N_{a,b,c}=\bigl\{(z_1&,z_2,z_3)\in\C^3:\ms{z_1}-\ms{z_2}=2a,
\;\> x\in\R,\;\> y\in(-R,R),\\
&z_1z_2=v_{a,b}(x,y)+iy,\quad z_3=x+iu_{a,b}(x,y)+ic\bigr\}/\Z.
\end{split}
\label{fi7eq8}
\e
Then $N_{a,b,c}$ is a noncompact SL 3-fold without boundary in $\C^3/\Z$,
which is nonsingular if $a\ne 0$, by Proposition~\ref{fi4prop1}.
\label{fi7def2}
\end{definition}

\begin{theorem} In the situation above, distinct\/ $N_{a,b,c}$
are disjoint. Define
\e
V=\bigl\{(z_1,z_2,z_3)\in\C^3:\md{\Im z_1z_2}<R\bigr\}/\Z.
\label{fi7eq9}
\e
Then there exists a continuous, surjective $F:V\ra\R^3$
with\/ $F^{-1}(a,b,c)=N_{a,b,c}$ for all\/ $(a,b,c)\in\R^3$.
Thus, $F$ is a special Lagrangian fibration of\/~$V$.
\label{fi7thm5}
\end{theorem}

\begin{proof} Clearly we have $N_{a,b,c}\subset V$ for each
$(a,b,c)\in\R^3$, where $N_{a,b,c}$ and $V$ are as in
\eq{fi7eq8} and \eq{fi7eq9}. Let $(z_1,z_2,z_3)\Z\in V$. We
shall show that there exists a unique $(a,b,c)\in\R^3$ such
that $(z_1,z_2,z_3)\Z\in N_{a,b,c}$. Let $x=\Re z_3$,
$y=\Im z_1z_2$ and $2a=\ms{z_1}-\ms{z_2}$. Consider the function
$b\mapsto v_{a,b}(x,y)$. By part (i) of Definition \ref{fi7def2}
and the last part of Theorem \ref{fi7thm3} this is a
{\it continuous} function.

If $b<b'$ then part (iv) of Definition \ref{fi7def2} gives
$\phi^\pm_{a,b}<\phi^\pm_{a,b'}$, and so Theorem \ref{fi7thm4}
gives $v_{a,b}<v_{a,b'}$ on $S$. Thus $b\mapsto v_{a,b}(x,y)$
is {\it strictly increasing}. And as the maximum and minimum
of $v_{a,b}$ on $S/\Z$ is achieved on $\pd(S/\Z)$ and so is a
value of $\phi^\pm_{a,b}$, part (v) of Definition \ref{fi7def2}
implies that $v_{a,b}(x,y)\ra\pm\iy$ as~$b\ra\pm\iy$.

Hence by the Intermediate Value Theorem there exists a unique
$b\in\R$ such that $v_{a,b}(x,y)=\Re z_1z_2$, and then $c$ is
given by $c=\Im z_3-u_{a,b}(x,y)$. It is easy to see from
\eq{fi7eq8} that $(z_1,z_2,z_3)\Z\in N_{a,b,c}$, and that this
is the only $(a,b,c)\in\R^3$ for which this holds. Therefore
distinct $N_{a,b,c}$ are disjoint, and $V=\bigcup_{(a,b,c)\in
\R^3}N_{a,b,c}$. The remainder of the proof follows
Theorem~\ref{fi4thm3}.
\end{proof}

Here is a simple, explicit example.

\begin{example} In the situation above, define $\phi^+_{a,b}(x)=
\phi^-_{a,b}(x)=b$ for all $a,b,x\in\R$. Then parts (i)--(v)
of Definition \ref{fi7def2} hold. It is easy to verify that
$v_{a,b}(x,y)\equiv b$ and $u_{a,b}(x,y)\equiv 0$. Hence
\begin{align*}
N_{a,b,c}=\bigl\{(z_1,z_2,z_3)\in\C^3:\,&\ms{z_1}-\ms{z_2}=2a,\quad
\Re z_1z_2=b,\quad
\Im z_3=c,\\
&\Im z_1z_2\in(-R,R)\bigr\}/\Z.
\end{align*}
It readily follows that the special Lagrangian fibration
$F:V\ra\C^3$ of Theorem \ref{fi7thm5} is given by
\e
F\bigl((z_1,z_2,z_3)\Z\bigr)=
\bigl(\ha\ms{z_1}-\ha\ms{z_2},\Re z_1z_2,\Im z_3\bigr).
\label{fi7eq10}
\e
Note that $F$ is a {\it smooth} SL fibration.

It is easy to show that $N_{a,b,c}$ is singular if and only if
$a=b=0$, so that the discriminant of $F$ is $\De=\bigl\{(0,0,c):
c\in\R\bigr\}$. This is of {\it codimension two} in $\R^3$, as
in Proposition \ref{fi3prop2}. The singular set of $N_{0,0,c}$
is $\bigl\{(0,0,x+ic):x\in\R\bigr\}/\Z$, which is a {\it circle}
${\cal S}^1$ in~$\C^3/\Z$.
\label{fi7ex1}
\end{example}

\subsection{A 1-parameter family of SL fibrations}
\label{fi73}

We now construct a 1-parameter family of SL fibrations $F^t$
for~$t\in[0,1]$.

\begin{example} Fix $P=2\pi$, let $t\in[0,1]$, and define
$\phi^+_{a,b}(x)=\phi^-_{a,b}(x)=b+t\cos x$ for all $a,b,x\in\R$.
Then parts (i)--(v) of Definition \ref{fi7def2} hold. Let $F^t$
be the SL fibration constructed in Theorem \ref{fi7thm5} using
this data. This gives a 1-{\it parameter family} of SL fibrations
$F^t:V\ra\C^3$ for $t\in[0,1]$, where $F^0$ coincides with the
smooth SL fibration of \eq{fi7eq10}. Use the notation $u_{a,b}^t,
v_{a,b}^t$ and $N_{a,b,c}^t$ in the obvious way. By the last part
of Theorem \ref{fi7thm3} the $N^t_{a,b,c}$ depend continuously on
$t$, and thus the $F^t$ depend continuously on~$t$.

From above, the singular fibres of $F^0$ are $N^0_{0,0,c}$, and
each is singular along a circle ${\cal S}^1$. However, for
$t>0$ things are different. By Theorem \ref{fi4thm5}, a
necessary condition for $N^t_{0,b,c}$ to have non-isolated
singularities is that $v_{0,b}^t(x,-y)\equiv -v_{0,b}^t(x,y)$.
But putting $y=R$ then gives $\phi^-_{0,b}(x)\equiv
-\phi^+_{0,b}(x)$, which is false. Thus when $t>0$, each
$N^t_{0,b,c}$ has isolated singularities, and only finitely
many of them by compactness.
\label{fi7ex2}
\end{example}

This will be our local model of how to deform a smooth SL
fibration to a non-smooth SL fibration. Here are some facts
we will need.

\begin{proposition} In Example \ref{fi7ex2}, for all\/ $t\in(0,1]$
and\/ $a,b\in\R$ we have
\begin{itemize}
\setlength{\itemsep}{0pt}
\setlength{\parsep}{0pt}
\item[{\rm(a)}] $x\mapsto v_{a,b}^t(x,y)$ has period\/ $2\pi$
and is strictly decreasing on $[0,\pi]$ and strictly increasing
on $[\pi,2\pi]$, for all\/~$y\in[-R,R]$.
\item[{\rm(b)}] $u^t_{a,b}(x,y)>0$ on $(0,\pi)\t(0,R]$ and\/
$(\pi,2\pi)\t[-R,0)$, and\/ $u^t_{a,b}(x,y)<0$ on $(0,\pi)
\t[-R,0)$ and\/~$(\pi,2\pi)\t(0,R]$.
\end{itemize}
\label{fi7prop1}
\end{proposition}

\begin{proof} By symmetries of the $\phi^\pm_{a,b}$ the
$v_{a,b}^t$ satisfy $v_{a,b}^t(-x,y)\equiv v_{a,b}^t(x,y)$ and
$v_{a,b}^t(2\pi-x,y)\equiv v_{a,b}^t(x,y)$. Hence
$\frac{\pd}{\pd x}v_{a,b}^t(0,y)\equiv\frac{\pd}{\pd x}v_{a,b}^t
(\pi,y)\equiv 0$ for $a\ne 0$. Also $\frac{\pd}{\pd x}v_{a,b}^t
(x,\pm R)\equiv -t\sin x$. Thus $\frac{\pd}{\pd x}v_{a,b}^t
(x,\pm R)<0$ for $x\in(0,\pi)$ and~$t>0$.

Now by \cite[\S 8.4]{Joyc2}, $\frac{\pd}{\pd x}v_{a,b}^t$ satisfies
a kind of {\it maximum principle}. Applying this on the rectangle
$[0,\pi]\t[-R,R]$ and using the fact that $\frac{\pd}{\pd x}v_{a,b}^t$
is zero on two sides of the rectangle and negative on the other
two, we find that $\frac{\pd}{\pd x}v_{a,b}^t<0$ on $(0,\pi)\t[-R,R]$
when $t>0$ and $a\ne 0$. Similarly, $\frac{\pd}{\pd x}v_{a,b}^t>0$ on
$(\pi,2\pi)\t[-R,R]$ when $t>0$ and $a\ne 0$. This implies part (a)
when~$a\ne 0$.

Also, by \eq{fi4eq5} we see that $\frac{\pd}{\pd y}u_{a,b}^t>0$ on
$(0,\pi)\t[-R,R]$ and $\frac{\pd}{\pd y}u_{a,b}^t<0$ on $(\pi,2\pi)
\t[-R,R]$. By symmetries of the $\phi^\pm_{a,b}$ the $u_{a,b}^t$
satisfy $u_{a,b}^t(x,-y)\equiv -u_{a,b}^t(x,y)$, and so
$u_{a,b}^t(x,0)\equiv 0$. Part (b) for $a\ne 0$ then follows by
integration on the line segment~$\{x\}\t[0,y]$.

Taking the limit $a\ra 0$ in (a) shows that $x\mapsto v_{0,b}^t(x,y)$
is decreasing on $[0,\pi]$ and increasing on $[\pi,2\pi]$. If it were
not {\it strictly} decreasing or increasing then it would be constant
on some subinterval $[\al,\be]$ of $[0,2\pi]$ with $\al<\be$. Then
$(u_{0,b}^t,v_{0,b}^t)$ is constant on $[\al,\be]\t\{0\}$, and
\cite[Th.~7.8]{Joyc4} implies that it is constant on $S$, a
contradiction. This completes part~(a).

Similarly, taking the limit $a\ra 0$ in (b) shows that
$u^t_{0,b}(x,y)\ge 0$ on $(0,\pi)\t(0,R]$ and $(\pi,2\pi)\t
[-R,0)$ and $u^t_{0,b}(x,y)\le 0$ on $(0,\pi)\t[-R,0)$ and
$(\pi,2\pi)\t(0,R]$. But if $u^t_{0,b}(x,y)=0$ at any interior
point of these regions, then we can derive a contradiction using
the {\it strong maximum principle} \cite[Th.~3.5]{GiTr}. This
completes part~(b).
\end{proof}

We can now prove an analogue of Proposition \ref{fi6prop} for
Example~\ref{fi7ex2}.

\begin{proposition} In Example \ref{fi7ex2} there are continuous
$\al,\be:[0,1]\ra\R$ with
\e
F^t\bigl((0,0,0)\Z\bigr)=\bigl(0,\al(t),0\bigr)
\quad\text{and}\quad
F^t\bigl((0,0,\pi)\Z\bigr)=\bigl(0,\be(t),0\bigr),
\label{fi7eq11}
\e
and\/ $\al(0)=\be(0)=0$. For all\/ $t\in(0,1]$ these have the
properties that
\begin{itemize}
\setlength{\itemsep}{0pt}
\setlength{\parsep}{0pt}
\item[{\rm(i)}] If\/ $b\notin\bigl[\al(t),\be(t)\bigr]$ then
$(u_{0,b}^t,v_{0,b}^t)$ has no singularities in~$S$.
\item[{\rm(ii)}] $(u^t_{0,\al(t)},v^t_{0,\al(t)})$ has singularities
of multiplicity $2$ and maximum type at\/ $(2n\pi,0)$ for $n\in\Z$,
and no other singularities.
\item[{\rm(iii)}] $(u^t_{0,\be(t)},v^t_{0,\be(t)})$ has singularities
of multiplicity $2$ and minimum type at\/ $(\pi+2n\pi,0)$ for $n\in\Z$,
and no other singularities.
\item[{\rm(iv)}] If\/ $b\in\bigl(\al(t),\be(t)\bigr)$ there exists
$x\in(0,\pi)$ such that\/ $(u_{0,b}^t,v_{0,b}^t)$ has singularities
of multiplicity $1$ and increasing type at\/ $(2n\pi-x,0)$ for
$n\in\Z$, singularities of multiplicity $1$ and decreasing type
at\/ $(2n\pi+x,0)$ for $n\in\Z$, and no other singularities.
\end{itemize}
\label{fi7prop2}
\end{proposition}

\begin{proof} Let $t\in[0,1]$. As in the proof of Theorem
\ref{fi7thm5}, the function $b\mapsto v^t_{0,b}(0,0)$ is
continuous, strictly increasing, and tends to $\pm\iy$ as
$b\ra\pm\iy$. Thus by the Intermediate Value Theorem there
exists a unique value $\al(t)$ such that $v^t_{0,\al(t)}(0,0)=0$.
Since $u^t_{0,\al(t)}(0,0)=0$ by definition, this means that
$(0,0,0)\Z\in N^t_{0,\al(t),0}$ by \eq{fi7eq8}, and so
$F^t\bigl((0,0,0)\Z\bigr)=\bigl(0,\al(t),0\bigr)$ by
definition of~$F^t$.

Hence $\al:[0,1]\ra\R$ exists, and is continuous as $F^t$
depends continuously on $t$. Similarly, considering the
function $b\mapsto v^t_{0,b}(\pi,0)$ we find a unique value
$\be(t)$ such that $v^t_{0,\be(t)}(\pi,0)=0$, and $\be:[0,1]
\ra\R$ exists and is continuous. Also $F^0$ is given in
\eq{fi7eq10}, so $\al(0)=\be(0)=0$ follows from~\eq{fi7eq11}.

Let $t\in(0,1]$, and consider the function $x\mapsto v^t_{0,b}(x,0)$.
By part (a) of Proposition \ref{fi7prop1} this has period $2\pi$ with
a maximum at 0 and a minimum at $\pi$. But $b\mapsto v^t_{0,b}(0,0)$
is strictly increasing and zero when $b=\al(t)$. Hence the
maximum of $x\mapsto v^t_{0,b}(x,0)$ is {\it negative} when
$b<\al(t)$, {\it zero} when $b=\al(t)$ and {\it positive} when
$b>\al(t)$. In particular, if $b<\al(t)$ then $v^t_{0,b}(x,0)<0$
for all $x$, so that $(u_{0,b}^t,v_{0,b}^t)$ has no singularities.
This proves half of part~(i).

By similar arguments using part (a) of Proposition \ref{fi7prop1}
and the Intermediate Value Theorem, we easily find that the
singularities of $(u_{0,b}^t,v_{0,b}^t)$, which are the zeroes
of $x\mapsto v^t_{0,b}(x,0)$, are as given in parts (i)--(iv),
and the {\it type} of each singularity also follows from part
(a) of Proposition \ref{fi7prop1}. Finally, to identify the
{\it multiplicity} of each singularity $(x,0)$ we can use part
(b) of Proposition \ref{fi7prop1}, as knowing the sign of
$u_{0,b}^t$ constrains the winding number of $(u_{0,b}^t,v_{0,b}^t)$
about 0 along $\ga_\ep(x,0)$ in Definition~\ref{fi4def3}.
\end{proof}

This gives the {\it discriminant} of~$F^t$.

\begin{corollary} In Example \ref{fi7ex2}, the discriminant of\/ $F^t$ is
\e
\De^t=\bigl\{(0,b,c):b\in [\al(t),\be(t)],\quad c\in\R\bigr\}\subset\R^3,
\label{fi7eq12}
\e
and the set of singular points is $\bigl\{(0,0,z_3):z_3\in\C\bigr\}/\Z$.
For $t\in(0,1]$, if\/ $b=\al(t)$ or $b=\be(t)$ then $N_{0,b,c}$ has
one singular point, and if\/ $\al(t)<b<\be(t)$ then $N_{0,b,c}$ has
two singular points.
\label{fi7cor}
\end{corollary}

Note that the set of singular points of $F^t$ in $\C^3$ is
{\it independent of}\/~$t$.

\subsection{Discussion}
\label{fi74}

Example \ref{fi7ex2} constructs a continuous 1-parameter family
of SL fibrations $F^t:V\ra\R^3$ for $t\in[0,1]$, where $F^0$ is
the {\it smooth} SL fibration given explicitly in \eq{fi7eq10},
but $F^t$ is not smooth for $t\in(0,1]$. Thus, this provides a
local model for how to continuously deform a smooth SL fibration
to a non-smooth SL fibration.

We have now seen several examples of non-smooth SL fibrations,
and a local mechanism for deforming smooth SL fibrations to
non-smooth ones. This justifies the following:

\begin{conjecture} Generic SL fibrations $f:M\ra B$ of almost
Calabi--Yau $3$-folds including fibres with singularities
are never smooth, but only piecewise smooth. Here we use
`generic' in the sense of Definition~\ref{fi3def4}.
\label{fi7conj}
\end{conjecture}

Suppose we have an almost Calabi--Yau 3-fold $M$ with a smooth
SL fibration $f:M\ra B$. What happens to the fibration if we
deform $M$ to a nearby generic almost Calabi--Yau 3-fold $\ti M$?
I believe that the SL fibration will still exist, at least on
most of $\ti M$, but that it will be only piecewise smooth.

For a smooth SL fibration the discriminant $\De$ is of codimension
two, and is expected to be a {\it graph}. In Example \ref{fi7ex2},
when $t=0$ and $F^0$ is smooth, the discriminant $\De^0$ is a
{\it line} in $\R^3$, of codimension two, but as $t$ increases
$\De^t$ thickens out continuously into a {\it ribbon} in $\R^3$,
of codimension one. In the same way, in a small generic deformation
$\ti M$ of $M$, I conjecture that the edges of the graph $\De$ in
$B$ thicken out into 2-dimensional ribbons in $\ti\De\subset B$. We
will discuss what happens near the vertices of $\De$ in~\S\ref{fi8}.

We can also relate the behaviour of $F^t$ to the local models of
\S\ref{fi5} and~\S\ref{fi6}:
\begin{itemize}
\setlength{\itemsep}{0pt}
\setlength{\parsep}{0pt}
\item When $\Re z_3\in(\pi,2\pi)$, the SL fibration $F^t$ for
$t\in(0,1]$ locally resembles the SL fibration $F$ of Theorem
\ref{fi5thm1} near $(0,0,z_3)\Z$, with singularities of
multiplicity 1 and increasing type.
\item When $\Re z_3\in(0,\pi)$, the SL fibration $F^t$ for
$t\in(0,1]$ locally resembles the SL fibration $F'$ of Theorem
\ref{fi5thm2} near $(0,0,z_3)\Z$, with singularities of
multiplicity 1 and decreasing type.
\item When $\Re z_3=0$, the SL fibration $F^t$ for $t\in(0,1]$
locally resembles the SL fibration $\hat F$ of Theorem
\ref{fi6thm2} near $(0,0,z_3)\Z$, with singularities of
multiplicity 2 and maximum type.
\item When $\Re z_3=\pi$, the SL fibration $F^t$ for $t\in(0,1]$
locally resembles an {\it analogue} of the SL fibration $\hat F$
of Theorem \ref{fi6thm2} near $(0,0,z_3)\Z$, with singularities
of multiplicity 2 and {\it minimum} type.
\end{itemize}
Thus, we have not found any new kinds of local singular behaviour
of SL fibrations in this section, we have just assembled those
already discussed in \S\ref{fi5} and \S\ref{fi6} in a model with
more interesting global topology.

\section{Global behaviour of SL fibrations}
\label{fi8}

We are now ready to state our picture (still conjectural and incomplete) 
of what SL fibrations of generic almost Calabi--Yau 3-folds look like,
if indeed they exist. Rather than starting from scratch, we begin in
\S\ref{fi81} by reviewing the elegant picture of smooth SL fibrations
$f:M\ra B$, which has been built up largely by Mark Gross and Wei-Dong
Ruan. Then in \S\ref{fi82} we explain how to modify the Gross--Ruan
picture under a small generic deformation of $M$. Finally, in
\S\ref{fi83} we draw some conclusions about the SYZ Conjecture.

\subsection{The Gross--Ruan picture of smooth SL fibrations}
\label{fi81}

Here is a review of the expected properties of smooth SL 
fibrations of Calabi--Yau 3-folds. Our principal sources are Gross 
\cite[\S 3]{Gros2} and Ruan \cite[\S 7]{Ruan2} for the topology 
of the singular fibres, and Gross \cite[\S 1]{Gros3} and Ruan 
\cite[\S 9]{Ruan3} for the monodromy matrices. A concise statement 
may be found in the `Precise SYZ mirror conjecture' of 
Ruan~\cite[\S 9]{Ruan3}. 

Let $f:M\ra{\cal S}^3$ be a smooth SL fibration, with 
fibres $N_b=f^{-1}(b)$, and generic fibre $T^3$. For generic such 
fibrations, the discriminant $\De$ is thought to be a {\it trivalent 
graph}, made up of smooth edges, and vertices of two kinds, which we 
shall refer to as {\it positive} and {\it negative}. The topology
and local monodromy for each kind of singular fibre are as follows.

\subsubsection*{(a) Edges}
Let $\ga$ be an edge in $\De$, and $b\in\ga$. Then $N_b$
has the topology of $T^3$ with $T^2$ collapsed to an ${\cal S}^1$,
and may be written $\Si\t{\cal S}^1$, where $\Si$ is a $T^2$ with an
${\cal S}^1$ collapsed to a point, or equivalently an ${\cal S}^2$ with
two points identified. These fibres are called type $(2,2)$ by Gross
and type $I$ by Ruan. They have Euler characteristic zero.

The monodromy about each edge $\ga$ in $\De$, acting on $H_1(T^3;\Z)
\cong\Z^3$, is 
\e
\begin{pmatrix} 1 & 1 & 0 \\ 0 & 1 & 0 \\ 0 & 0 & 1 \end{pmatrix}
\label{fi8eq1}
\e
with respect to a suitable basis of~$H_1(T^3;\Z)$. 

\subsubsection*{(b) Positive vertices}
Let $b$ be a positive vertex in $\De$. Then $N_b$
has the topology of $T^3$ with $T^2$ collapsed to a point. It
has Euler characteristic 1. These fibres are called type $(1,2)$ 
by Gross and type $III$ by Ruan.

The monodromies around the three edges $\ga_1,\ga_2,\ga_3$ meeting 
at $b$ are
\e
\begin{pmatrix} 1 & 0 & 0 \\ 1 & 1 & 0 \\ 0 & 0 & 1 \end{pmatrix},\quad
\begin{pmatrix} 1 & 0 & 0 \\ 0 & 1 & 0 \\ -1 & 0 & 1 \end{pmatrix}
\quad\text{and}\quad
\begin{pmatrix} 1 & 0 & 0 \\ -1 & 1 & 0 \\ 1 & 0 & 1 \end{pmatrix}
\label{fi8eq2}
\e
with respect to a suitable basis of~$H_1(T^3;\Z)$. 

The smooth SL fibration of Example \ref{fi4ex1} is 
a local model for the fibration $f$ near the singular point of a
positive singular fibre.

\subsubsection*{(c) Negative vertices}
Let $b$ be a negative vertex in $\De$. Then Ruan
\cite[\S 7]{Ruan2} gives two different possible topologies for 
$N_b$, which he calls type $II$ and type $\ti{II}$. His
type $\ti{II}$ topology agrees with Gross' proposed type 
(2,1) fibre~\cite[\S 3]{Gros2}. 

Both fibres are constructed by taking a fibration $\pi:T^3\ra T^2$ 
with fibre ${\cal S}^1$, and collapsing the fibres to points over a 
graph $\Ga$ in $T^2$. In the type $II$ case $\Ga$ has three edges and 
two vertices, and in the type $\ti{II}$ case it has two edges 
and one vertex. In both cases $N_b$ has Euler characteristic~$-1$. 

The monodromies around the three edges $\ga_1,\ga_2,\ga_3$ meeting 
at $b$ are
\e
\begin{pmatrix} 1 & 1 & 0 \\ 0 & 1 & 0 \\ 0 & 0 & 1 \end{pmatrix},\quad
\begin{pmatrix} 1 & 0 & -1 \\ 0 & 1 & 0 \\ 0 & 0 & 1 \end{pmatrix}
\quad\text{and}\quad
\begin{pmatrix} 1 & -1 & 1 \\ 0 & 1 & 0 \\ 0 & 0 & 1 \end{pmatrix},
\label{fi8eq3}
\e
with respect to a suitable basis of~$H_1(T^3;\Z)$. 

At present, to the author's knowledge, there is no known local model 
for a smooth special Lagrangian fibration (or even a smooth Lagrangian
fibration) in the neighbourhood of a codimension three singular point
of a negative singular fibre. The author conjectures that no such local
model exists. If this is the case then smooth SL
fibrations may not exist on general Calabi--Yau 3-folds, even with a
very nongeneric choice of almost Calabi--Yau metric.
\medskip

We will refer to the singular fibres over positive and negative vertices 
as {\it positive} and {\it negative singular fibres} respectively. Our 
notation of positive and negative vertices was suggested by David 
Morrison, and refers to the sign of the Euler characteristic of the
singular fibres. Gross' notation refers to the Betti numbers $(b^1,b^2)$ 
of the singular fibres.

\subsection{Modification of this picture for generic ACY 3-folds}
\label{fi82}

Now we shall something about what special Lagrangian fibrations
of generic almost Calabi--Yau 3-folds might look like. Suppose
we start with a smooth Gross--Ruan fibration $f:M\ra B$, either
of a nongeneric almost Calabi--Yau 3-fold or of the degenerate
large complex structure limit, and make a small perturbation to a
generic almost Calabi--Yau 3-fold. What happens to the fibration?

Near a nonsingular fibre $N_b=f^{-1}(b)$ of $f$, the fibration should 
remain nonsingular, and the local geometry unchanged. The interesting 
question is what happens to the singular fibres of $f$. The following 
is the author's best guess, on the assumption that SL 
fibrations are well-behaved in the generic case. We preface it with
some remarks on monodromy and coordinates on the moduli space.

Let $f:M\ra B$ be an SL fibration. By Theorem \ref{fi2thm1}, near a
nonsingular fibre $N_b\cong T^3$ the moduli space of deformations of
$N_b$ is isomorphic to $H^1(N_b;\R)\cong\R^3$. But this moduli space
is $B$, and so near any point in $B\sm\De_f$ we have natural affine
coordinates modelled on~$H^1(T^3;\R)$.

However, near a singular fibre $N_b$ the situation is more complicated
because of the monodromy action. Let $N_{b'}$ be a nonsingular fibre
near $N_b$. Let $\Ga_b$ be the set of monodromies of loops in $B\sm\De_f$
based at $b'$ and staying in a small neighbourhood of $b$. Then $\Ga_b$
is a group acting on $H_1(N_{b'};\Z)$ and $H^1(N_{b'};\R)$. Roughly 
speaking, near $b$ we can regard $B$ as a kind of quotient of 
$H^1(N_{b'};\R)$ by $\Ga_b$, so that $B$ is a kind of orbifold, with 
the topology of a 3-manifold, but not the smooth structure. 

In what follows, as long as we make use of only $\Ga_b$-invariant
objects, we can think of $B$ as being locally like $\R^3$ and 
mostly ignore the monodromy action. We shall represent elements of 
$H_1(N_{b'};\Z)$ by column vectors, and elements of $H^1(N_{b'};\R)$ 
by row vectors, upon which the monodromy matrices of equations 
\eq{fi8eq1}--\eq{fi8eq3} act by left and right multiplication 
respectively. 

Here is a conjectural picture of how generic SL fibrations work
around perturbations of the Gross--Ruan singular fibres described
in (a)--(c) of~\S\ref{fi81}.

\subsubsection*{(a) Edges}
The author conjectures that under small deformations, the `edges'
$\ga$ in the Gross--Ruan picture will thicken out into thin `ribbons'
$R$ of the kind described in \S\ref{fi7}. They are closed subsets of 
hyperplanes in $B$ defined locally by $[\om]\cdot[D]=0$, where $[\om]$ 
is the relative de Rham cohomology class in $H^1(M,N_b;\R)$ and $[D]$ a
relative homology class in $H_1(M,N_b;\Z)$ depending on the edge, which
will be represented by an even number of holomorphic discs $D$ for
generic $b\in B$, as we discussed in \S\ref{fi51} and~\S\ref{fi61}.

Choose an identification $H_1(N_b;\Z)\cong\Z^3$ such that the monodromy
around $\ga$ is as in \eq{fi8eq1}. Then calculation using the local
model of \S\ref{fi7} shows that $[\pd D]\in H_1(N_b;\Z)$ should be
identified with $\pm(1\,0\,0)^T$ in $\Z^3$. Also, $B$ is locally
identified with $H^1(N_b;\R)\cong\R^3$ up to monodromy, and $R$ lies
in the monodromy-invariant hyperplane~$\bigl\{(0,x_2,x_3):x_j\in
\R\bigr\}$.

The situation described in \S\ref{fi7}, in which the generic singular
fibre has two singular points, is only the simplest possibility. In
general we expect the generic singular fibre to contain an even number 
of singular points, divided equally into two kinds. In codimension one 
on the ribbon these singular points can appear or disappear in pairs of
different kinds, and the edge of the ribbon is where the last two singular 
points disappear.

\subsubsection*{(b) Positive vertices}
For positive vertices in the Gross--Ruan picture, the
monodromy matrices of \eq{fi8eq2} all fix the vectors
\begin{equation*}
{\bf v}_1=\begin{pmatrix} 0 \\ 1 \\ 0\end{pmatrix}, \quad
{\bf v}_2=\begin{pmatrix} 0 \\ 0 \\ -1\end{pmatrix} 
\quad\text{and}\quad
{\bf v}_3=\begin{pmatrix} 0 \\ -1 \\ 1\end{pmatrix}
\end{equation*}
in $H_1(N_{b'};\Z)$ and the direction $(1\,0\,0)$ in~$H^1(N_{b'};\R)$. 

\begin{figure}[htb]
\centerline{$\splinetolerance{.7pt}
\begin{xy}
0;<1mm,0mm>:
,(0,0)*{\bullet}
,(0,-7.5)*{\bullet}
,(0,-15)*{\bullet}
,(36,12.5)*!L{R_1}
,(-16,-32.5)*!R{R_2}
,(-1,-7.5)*!R{\scriptstyle b_0}
,(-26,17.5)*!R{R_3}
,(0,0);(0,-15)**@{-}
,(0,0);(0,-15)**\crv{(10,-10)&(10,-25)}
?(.06)="a"
?(.13)="b"
?(.20)="c"
?(.27)="d"
?(.35)="e"
?(.26)="p"
?(.43)="q"
?(.57)="r"
?(.73)="s"
,(0,0);(7,4)**\crv{(6,10)&(6.5,6)}
?(.7)="x"
,(7,4);(0,-15)**\crv{~*=<3pt>{.}(7.5,2)&(6,-5)}
,(0,0);(-9,-15)**\crv{(-7,-4)&(-11,-9)}
?(.1)="t"
?(.2)="u"
?(.3)="v"
?(.4)="w"
?(.5)="y"
?(.58)="z"
?(.77)="zz"
?(.92)="zzz"
,(-9,-15);(0,-15)**\crv{~*=<3pt>{.}(-7,-21)&(-3.5,-17)}
,(0,0);(35,20)**@{-}
,"a";(35,17)**@{.}
,"b";(35,14)**@{.}
,"c";(35,11)**@{.}
,"d";(35,8)**@{.}
,"e";(35,5)**@{-}
,(0,-3);"p"**@{.}
,(0,-6);"q"**@{.}
,(0,-9);"r"**@{.}
,(0,-12);"s"**@{.}
,(3.5,2);"x"**@{.}
,(0,-15);"e"**\crv{~*=<3pt>{}~**{.}}
,(0,0);(-15,-25)**@{-}
,(0,-3);(-15,-28)**@{.}
,(0,-6);(-15,-31)**@{.}
,(0,-9);(-15,-34)**@{.}
,(0,-12);(-15,-37)**@{.}
,(0,-15);(-15,-40)**@{-}
,(0,0);(-25,25)**@{-}
,"t";(-25,22)**@{.}
,"u";(-25,19)**@{.}
,"v";(-25,16)**@{.}
,"w";(-25,13)**@{.}
,"y";(-25,10)**@{-}
,"y";(0,-15)**\crv{~*=<3pt>{.}(-6,-9)}
,"z";(-2.4,-4)**@{.}
,"zz";(-4.8,-8)**@{.}
,"zzz";(-7.2,-12)**@{.}
\end{xy}$}
\caption{Discriminant locus near a perturbation of a positive vertex}
\label{fi8fig1}
\end{figure}
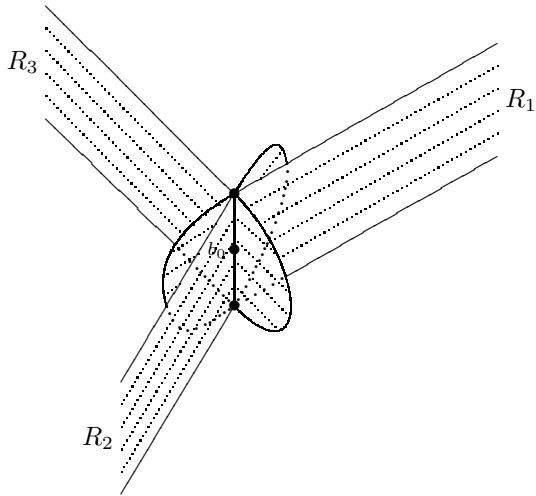

In a generic perturbation of a Gross--Ruan fibration near a positive
vertex, the three edges in $\De_f$ should thicken out into `ribbons' 
$R_1,R_2,R_3$ lying in the three hyperplanes
\begin{gather*}
H_1=\bigl\{(x_1,0,x_3):x_j\in\R\bigr\},\qquad
H_2=\bigl\{(x_1,x_2,0):x_j\in\R\bigr\}\\
\text{and}\qquad H_3=\bigl\{(x_1,x_2,x_3):x_j\in\R,\quad
x_2=x_3\bigr\},
\end{gather*}
which are the hyperplanes dual to ${\bf v}_1,{\bf v}_2,{\bf v}_3$, 
and intersect in $\bigl\{(x_1,0,0):x_1\in\R\bigr\}$. The ribbons
$R_1,R_2,R_3$ intersect in a bounded subinterval of this line, as
sketched in Figure~\ref{fi8fig1}.

There are two obvious ways for this to happen, in which either
$R_1\cap R_2\cap R_3$ is part of the boundary of each $R_j$, or
the ribbons $R_j$ extend a little way beyond their intersection
$R_1\cap R_2\cap R_3$. The author thinks that the latter option
is what actually happens, as in Figure~\ref{fi8fig1}.

For generic points in the intersection $R_1\cap R_2\cap R_3$ the
singularities of the fibres are just finitely many points modelled 
locally on the $T^2$-cone $C$ of \eq{fi5eq12}. These are divided
into three kinds, corresponding to the ribbons $R_1,R_2,R_3$,
according to the homology class of the ${\cal S}^1$ in $T^3$
that collapses to a point.

However, at certain special points $b_0$ in $R_1\cap R_2\cap R_3$ 
there will be a new kind of codimension three singularity, when two 
or three of these singular points of different kinds come together.
The author does not have a local model for this singularity, but
topologically it may involve a cone on a genus 2 surface. There 
must be at least one such singular fibre, as it is necessary for
the monodromy to work out.

When $N_{b'}$ is a generic nonsingular fibre near the ribbon $R_j$,
there should exist an even number of holomorphic discs $D_j$ in $M$
whose boundary $\pd D_j$ in $N_{b'}$ has homology class $\pm{\bf v}_j$
in $H_1(N_{b'};\Z)\cong\Z^3$. Singularities develop when the area of
$D_j$ shrinks to zero, which happens on the hyperplane $H_j$ in~$B$.

\subsubsection*{(c) Negative vertices}
\begin{figure}[htb]
\centerline{$\splinetolerance{.8pt}
\begin{xy}
0;<1mm,0mm>:
,(35,20);(-25,25)**\crv{(4,4)&(-4,4)}
?(.95)="a"
?(.85)="b"
?(.75)="c"
?(.65)="d"
?(.55)="e"
?(.45)="f"
?(.35)="g"
?(.25)="h"
?(.15)="i"
?(.05)="j"
,(35,5);(-15,-40)**\crv{(8,0)&(-8,-8)}
?(.08)="k"
?(.2)="l"
?(.33)="m"
?(.5)="n"
?(.8)="o"
,(-21,-28);(-31,13)**\crv{(-7,-4)&(-7,1)}
?(.94)="p"
?(.84)="q"
?(.73)="r"
?(.58)="s"
?(.05)="t"
,"a";"p"**@{.}
,"b";"q"**@{.}
,"c";"r"**@{.}
,"d";"s"**@{.}
,"e";"t"**@{.}
,"j";"k"**@{.}
,"i";"l"**@{.}
,"h";"m"**@{.}
,"g";"n"**@{.}
,"f";"o"**@{.}
,(35,20);(35,5)**@{--}
,(-25,25);(-31,13)**@{--}
,(-21,-28);(-15,-40)**@{--}
\end{xy}$}
\caption{Discriminant locus near a perturbation of a negative vertex}
\label{fi8fig2}
\end{figure}
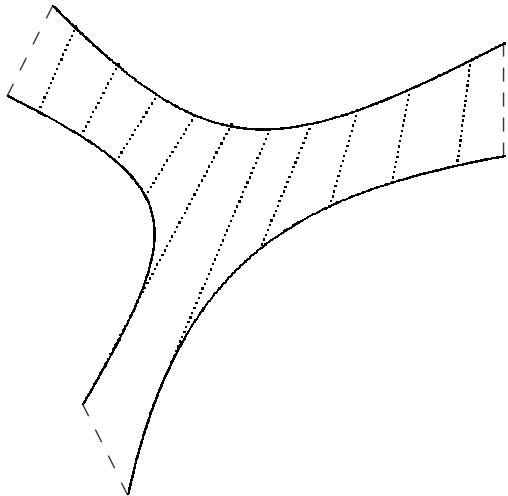

For negative vertices, the monodromy matrices of \eq{fi8eq3}
all fix the vector $(1\,0\,0)^T$ in $H_1(N_{b'};\Z)$ and the hyperplane 
$\bigl\{(0,x_2,x_3)\!:\!x_j\!\in\!\R\bigr\}$ in $H^1(N_{b'};\R)$. 
In a generic perturbation of a Gross--Ruan fibration near a negative
vertex, the three edges in $\De_f$ should thicken out into `ribbons' 
which all lie in the same hyperplane $H$ in $B$, isomorphic to 
$\bigl\{(0,x_2,x_3):x_j\in\R\bigr\}$ in $H^1(N_{b'};\R)$. The three 
ribbons merge together to make a letter $Y$ shape in $H$, as sketched 
in Figure~\ref{fi8fig2}. 

When $N_{b'}$ is a generic nonsingular fibre near this part of $\De_f$,
there should exist an even number of homologous holomorphic discs $D$
in $M$ whose boundaries $\pd D$ in $N_{b'}$ have homology class
$\pm(1\,0\,0)^T$ in $H_1(N_{b'};\Z)\cong\Z^3$. Singularities develop
when the area of $D$ shrinks to zero, which happens on the hyperplane
$H$ in~$B$.

Calculations by the author, along the lines of \S\ref{fi7} but 
more complicated, show that one can put together a fibration with 
the topological properties we want using only the local models of 
\S\ref{fi5} and \S\ref{fi6}. There is no need to include any other
kind of singular point.
\medskip

The author is fairly confident about parts (a) and (c), but less happy
about part (b). In fact, Ruan's Lagrangian fibrations by gradient flow
look quite like parts (a) and (c) in the relevant regions. Another
option in part (b) is that there could be a new kind of codimension 
two singularity along the line segment~$R_1\cap R_2\cap R_3$.

\subsection{Conclusions}
\label{fi83}

If the speculations of \S\ref{fi82} are correct, they have important 
consequences for the SYZ Conjecture. Positive and negative singular 
fibres are expected to be dual to one another under the mirror 
transform. That is, if we have dual smooth SL fibrations $f:M\ra B$ and 
$\hat f:\hat M\ra B$ as in the SYZ conjecture, then positive vertices 
in the discriminant $\De_f$ of $f$ in $B$ should coincide with negative
vertices in the discriminant $\De_{\smash{\hat f}}$ of $\hat f$, and 
vice versa. This follows as the monodromy matrices in \eq{fi8eq2} are
the transposes of those in~\eq{fi8eq3}.

However, after a small generic perturbation of $f$ and $\hat f$ near such 
a vertex in $B$, it is clear from Figures \ref{fi8fig1} and \ref{fi8fig2} 
that the discriminant loci $\De_f$ and $\De_{\smash{\hat f}}$ can no 
longer be identified, because they are not homeomorphic. On this basis 
we make the following conjecture.

\begin{conjecture} Let\/ $M,\hat M$ be generic mirror Calabi--Yau $3$-folds.
Then even if there do exist special Lagrangian fibrations $f:M\ra B$ 
and\/ $\hat f:\hat M\ra\hat B$, it is not in general possible to 
homeomorphically identify the bases $B$ and\/ $\hat B$ of the 
fibrations in a way that identifies the discriminants $\De_f$, 
$\De_{\smash{\hat f}}$ of\/ $f,\hat f$, and so that the nonsingular 
fibres of\/ $f,\hat f$ are $3$-tori with dual homology.
\label{fi8conj}
\end{conjecture}

This contradicts the version of the SYZ Conjecture stated in
the introduction, and some of the stronger forms of the SYZ
Conjecture that people have written down so far. If it is true
then it will limit the scope of any eventual final formulation
of the SYZ Conjecture.

My feeling is that while the SYZ Conjecture is clearly morally
true, it is probably not literally true of genuine special
Lagrangian fibrations of holonomy $\SU(3)$ Calabi--Yau 3-folds,
except in some limiting sense in the large complex structure limit.

Furthermore, I believe that the Gross--Ruan picture of smooth SL
fibrations is probably {\it asymptotically true} of general SL
fibrations in the large complex structure limit, so that in a
family of Calabi--Yau 3-folds with SL fibrations approaching the
large complex structure limit, the 2-dimensional discriminants
will collapse down onto 1-dimensional trivalent graphs.

Therefore, a better way to formulate the SYZ Conjecture might be
in terms of SL fibrations of 1-{\it parameter families} of mirror
Calabi--Yau 3-folds $M_t,\hat M_t$ for $t\in(0,\ep)$, which
both approach the large complex structure limit as $t\ra 0$.
A similar conclusion is reached by Gross in~\cite[\S 4]{Gros4}.

\end{document}